\documentclass[11pt]{amsart}
\usepackage{graphicx}
\usepackage{amssymb}
\usepackage{epstopdf}
\usepackage{amsmath}
\usepackage{amscd}
\usepackage{amsthm,amsfonts,amsxtra,amssymb,amscd}
\usepackage[all]{xy}
\usepackage{enumerate}
\usepackage{hyperref}

\newtheorem{definition}{Definition}
\newtheorem{theorem}{Theorem}[section]

\newtheorem{lemma}[theorem]{Lemma}
\newtheorem{proposition}[theorem]{Proposition}
\newtheorem{corollary}[theorem]{Corollary}

\newtheorem{remark}[theorem]{Remark}
\newtheorem{remarks}[theorem]{Remarks}

 \title{On the geometry of toric arrangements}

\author{C. De Concini}\address{Dip. Mat. Castelnuovo, Univ. di Roma La
Sapienza, Rome, Italy}\email{deconcin@mat.uniroma1.it}
\author{C. Procesi}\address{Dip. Mat. Castelnuovo, Univ. di Roma La
Sapienza, Rome, Italy}\email{procesi@mat.uniroma1.it}
\thanks{The authors are partially supported by the Cofin 40
\%, MIUR}
\begin{document}

\begin{abstract} Motivated by  the counting formulas of integral polytopes  as in  Brion--Vergne \cite{BV}, \cite{BV1} and Szenes--Vergne  \cite{SV1}, we start to lie the foundations of a theory for toric arrangements, which may be considered as the {\it periodic version} of the theory of hyperplane arrangements.\end{abstract}

\dedicatory{Dedicated to Vladimir Drinfeld on the occasion of his 50th birthday}\maketitle

  \section{Introduction}
  \subsection{The knapsack problem}
One of the aims of this paper is to understand in a more algebraic geometric language, the counting formulas of \cite{BV}, \cite{BV1},  \cite{SV1} for integral points in polytopes.

We are given a family  $\Delta$ of vectors in a lattice $\Lambda$  and, for $\beta\in\Lambda$ we want to count the number  $c_\beta$  of solutions of the equation:
$$\sum_{\alpha\in\Delta}n_\alpha\alpha=\beta,\ n_\alpha\in\mathbb N.$$
This is clearly the value of a {\it partition function}, or in other words, the number of integral points of the convex polytope
$$\Pi_\beta:=\{(x_\alpha)\,| \sum_{\alpha\in\Delta}x_\alpha\alpha=\beta,\ 0\leq x_\alpha\in\mathbb R\}.$$
One usually restricts to the case in which $\Delta$ is all on one side of a hyperplane so that the polytope is bounded and we have finitely many solutions. In this case 
the numbers $c_\beta$  are finite and are the coefficients of the generating series:
\begin{equation}
f_\Delta:=\sum_\beta c_\beta e^\beta=\prod_{\alpha\in \Delta}\frac{1}{1-e^\alpha}\end{equation}

Thus the problem is to understand an {\it inversion formula}, the formula  which computes these numbers $c_\beta$ from $f_\Delta$.\smallskip

A prototype example is the {\it knapsack problem}, i.e. the same problem where now the elements $\alpha$ are just positive integers  $h_i$, (cf. the paper of E.T. Bell \cite{Be}). In this setting,  we ask, given an integer $n$ to compute the coefficient of $x^n$ in the function $\prod_i{1\over 1-x^{h_i}}$  or to compute the residue  ${1\over 2\pi i}\oint \prod_i{x^{-n-1}\over 1-x^{h_i}}dx$ over a small circle.   The function $\prod_i{x^{-n-1}\over 1-x^{h_i}}$ has poles in 0 and in the $m-$th roots of 1, where $m$ is the least common multiple of the numbers $h_i$. In other words, if $\zeta=e^{2\pi i/m}$  one can write $\prod_i{1\over 1-x^{h_i}}=\prod_{i=1}^m{1\over (1-\zeta^ix)^{b_i}}$  where the integers $b_i$ are easily computed from the $h_j$. By standard residue theory we have:
\begin{equation}\label{knap} {1\over 2\pi i}\oint \prod_i{x^{-n-1}\over 1-x^{h_i}}dx=-\sum_{j=1}^m{1\over 2\pi i}\oint_{C_j}\prod_{i=1}^m{x^{-n-1}\over (1-\zeta^ix)^{b_i}}dx
\end{equation} where $C_j$ is a small circle around $\zeta^{-j}$. 

We are using the fact that the residue at $\infty$ is 0, due to the hypothesis made that $n>0$.\smallskip

For the
 term ${1\over 2\pi i}\oint_{C_j}\prod_{i=1}^m{x^{-n-1}\over (1-\zeta^ix)^{b_i}}dx
$ one makes a coordinate change $x=w+\zeta^{-j}$  and obtains $${1\over 2\pi i}\oint_{C_j}\prod_{i=1}^m{(w+\zeta^{-j})^{-n-1}\over (1-\zeta^{i-j}-\zeta^iw)^{b_i}}dw.
$$
Now  $\prod_{i=1,\ i\neq j}^m{1\over (1-\zeta^{i-j}-\zeta^iw)^{b_i}}$ is holomorphic around 0 and so it can be expanded as a power series $\sum_{h=0}^\infty a_{j,h}w^h$ whose terms can be computed, while $$(w+\zeta^{-j})^{-n-1}= \zeta^{j(n+1)}\sum_{k=0}^\infty (-1)^k\binom{n+k}{k}(\zeta^jw)^k.$$ Finally we have for the $j-$th term of (\ref{knap}):
$$-{(-1)^{b_j}\over 2\pi i}\oint_{C_j} \zeta^{j(n+1-b_j)}(\sum_{k=0}^\infty (-1)^k\binom{n+k}{k}\zeta^{jk}w^k)(\sum_{h=0}^\infty a_{j,h}w^h)w^{-b_j}dw= $$ $$-(-1)^{b_j}\zeta^{j(n+1-b_j)}\sum_{k+h=b_j-1}(-1)^k\zeta^{jk}\binom{n+k}{k}a_{j,h}.$$
 This formula, summed over all $j$, answers the question and exhibits the partition function requested as a sum of functions which are polynomials on the cosets modulo $m$. One calls such a function a {\it periodic polynomial} or {\it quasipolynomial}.
 
 In this paper we show that the counting formulas of \cite{BV}, \cite{BV1},  \cite{SV1}, which are the higher dimensional case of the knapsack problem,  can be interpreted in a similar way as sums of residues at the {\it poles}  of the generating function.  In order to develop this higher dimensional case, we shall study a trigonometric   analog  of the theory of hyperplane arrangements, which we call a {\it toric arrangement}, in which the main object of study is the open set of a torus complement of a finite set of cosets of  codimension 1 tori.  Versions  of these arrangements also known as toral arrangements have been introduced and studied in \cite {dl},  \cite {dou}, \cite{Le1}, \cite{Le2}, \cite{lo},  \cite {ma1},  \cite {ma2}, in particular in  \cite {lo},    Looijenga proves, by a completely different argument, using the degeneration of a spectral sequence, that the cohomology of this space has a natural filtration whose graded spaces coincide with the ones we describe in Theorem 4.2,  in particular the Betti numbers are computed. As for the algebra structure, contrary to what claimed in \cite{lo}, 2.4.3, the algebra is generated by the logarithmic  differentials if and only if the set of characters is {\it unimodular}, see \S  5.2 and  (\ref{unim}). Toric arrangements are   the main topic of this paper.
 
 There is a more analytic approach  to counting formulas, via splines. It is due to Dahmen, Micchelli cf. \cite{DM}, and it would be interesting to compare the two approaches.
 
 Finally there is some overlap, at least in the development of expansion formulas, with the work of  A. Szenes, \cite{S1}.

 \subsection{Toric arrangements} 
 
Let $T=hom_{\mathbb Z}(\Lambda,\mathbb C^*)$ be an algebraic torus over $\mathbb C$, of dimension $r$. $\Lambda$ its character group,   a rank $r$ lattice,
$A:=\mathbb C[T]=\mathbb C[\Lambda]$ its coordinate ring.
 \smallskip

 We consider, as analogue of a hyperplane, a translate of the kernel of a character $\chi$, i.e. the hypersurface
$H_{a,\chi}$ of equation $1-a\chi=0$ for a suitable non zero $a\in \mathbb C^*$. \smallskip

Given a finite subset  $\Delta\subset \mathbb C^*\times\Lambda$, we shall study the {\bf components of the arrangement $\mathcal R_\Delta$} generated by the
hypersurfaces $H_{a,\chi}$ as $(a,\chi)\in\Delta$. By this we shall mean {\it the set of all connected components} of all intersections of these
hypersurfaces. \smallskip

The complement of the union of these hypersurfaces will be denoted by $\mathcal A_\Delta$, it is an affine variety of coordinate ring $A[d^{-1}
],\ d:=\prod_{(a,\chi)\in\Delta}(1-a\chi).$\smallskip

  Our first concern will be to compute the cohomology of $\mathcal A_\Delta$. For this, let us remark
that  on the entire torus $T$ we have the closed forms $d\log\chi, \chi\in\Lambda$, spanning the vector space of invariant 1-forms isomorphic to $\Lambda\otimes\mathbb R$ and generating the cohomology of $T$. On the open set   $\mathcal A_{\Delta}$ 
we have then also the closed one forms  $d\log(1- a\chi),\ (a,\chi)\in\Delta $.     In general, unless the arrangement is {\it unimodular} (see \ref{unim}),    these two sets of 1-forms are not enough to generate cohomology (this corrects a statement in \cite{lo} 2.4.3).

  Nevertheless one can give a
satisfactory description, in term of explicit differential forms of  $H^*(\mathcal A_{\Delta},\mathbb C)$. One obtains then:

{\bf Theorem \ref{decomp}.} {\it
    For each integer $i\geq 0$, we have a (non canonical) decomposition, as $W$ runs over the components of the arrangement:}
    $$H^i(\mathcal A_{\Delta})=\oplus_{W} H^{i-{\rm codim }W}(W)
    \otimes  V_{W}.$$
  Each $W$ is isomorphic to a torus, so its cohomology is an exterior algebra on $r$ explicit forms, and   $V_{W}$ is a finite dimensional vector space, which depends on the combinatorics of the set $\Delta$ via the theory on non broken circuits and can be identified to  the top cohomology of the hyperplane arrangement defined by the differentials, at  a given point of $W$, of the functions $1-a\chi$ vanishing on $W$. \medskip

The result is particularly striking in the unimodular case, where we shall prove {\it formality}, i.e.

{\bf Theorem \ref{formality}}   {\it The subalgebra of the algebra of
differential forms generated by the 1-forms $d\log(1- a\chi), d\log \chi$ is isomorphic to the De Rham cohomology of $\mathcal A_\Delta$. }\smallskip

  The algebraic relations
between these forms resemble, but are more complicated, the relations of Orlik-Solomon (cf. 5.1) for the hyperplane case.  The theory of non broken circuits furnishes a linear
basis. \smallskip

The proof of these results will be based on a careful analysis of $\mathcal A_\Delta$. We shall use some simple ideas from $D-$modules to do a
precise bookkeeping.
\smallskip

The  reader unfamiliar with $D-$modules can find, in  the introductory book \cite{Co}  of S. C. 
Coutinho,  essentially all that we need.

   We have limited the analysis to De Rham cohomology. It is likely that a finer geometric analysis will produce a similar result (after normalizing the classes), for integral cohomology. 

\subsection{Residues}  The    second task is to define correctly the residues which will be integrals over suitable cycles. 
 Since we are in dimension $>1$ we are faced with the problem that the poles are on divisors with a complicated intersection pattern, this implies that we need to  use a model where the divisor has normal crossings as in \cite{dp}.  This is done in \S 6.

 Finally we have to understand which cycles contribute to the counting formula.  This is done in \S 7, it   is just a reinterpretation of results of Szenes--Vergne \cite{SV1} and it is connected to the  combinatorial problem of understanding the regions (chambers)  over which the partition function is a periodic polynomial. The answer is expressed by the notion of the Jeffrey--Kirwan residue.\smallskip

 With this analysis the counting formulas are explained  in a very similar way to the 1-dimensional case.

\vskip10pt
  \section{Periodic Weyl algebra}
\subsection{Differential operators}
    Given a torus $T$ with character group $\Lambda$, we shall now consider  the ring $D_T$ of differential operators over $T$. $D_T$ is generated by $A=\mathbb C[T]=\mathbb C[\Lambda]$ and by the invariant derivations $D_{\phi}$,
    $\phi\in\Lambda^*$ which are defined by $D_\phi(\chi):=\langle\phi,\chi\rangle\chi$, for any $\chi\in\Lambda$.  In terms of a basis of
    characters $\{\xi_1,\ldots ,\xi_m\}$, we have the basis of invariant derivations $\partial_i:=\xi_i{\partial\over\partial \xi_i}$.

    We shall think of  $D_T$ as a periodic version of  the Weyl algebra.  In explicit coordinates the
    {\it periodic } Weyl algebra is
$$W^p(r):=\mathbb C[\xi_1^{\pm 1},\dots,\xi_r^{\pm 1},{\partial\over\partial \xi_1},\dots, {\partial\over\partial \xi_r}]=\mathbb C[\xi_1^{\pm 1},\dots,\xi_r^{\pm 1},{\partial_1},\dots, {\partial_r}]$$

    If $\Lambda$ is the character group of  $T$ we shall also denote  $D_T=W(\Lambda)$.

     We filter the algebra $D_T$   by writing its elements as $$u:=\sum
 a_{i_1,\dots,i_r}\partial_1^{i_1}\dots\partial_r^{i_r},\quad  a_{i_1,\dots,i_r}\in A $$ and putting such an element in degree $\max\sum_{j=1}^r i_j$
 for which the coefficient $a_{i_1,\dots,i_r}$ is non 0.   The associated graded algebra is a polynomial ring $A[\eta_1,\dots,\eta_r]$ which can
 be viewed as the coordinate ring of the cotangent bundle to the torus $T$ (trivialized by the invariant forms).

   When we have a finitely generated  module $M$  one can define its {\it characteristic variety}, which is an important geometric invariant of the module,  applying to $D_T$ the same method used for modules over the Weyl algebra (cf. \cite{Co}). We define a {\it
   Bernstein filtration}, by choosing a finite set of generators $m_i$ and setting $M_j:=\sum_i(D_T)_jm_i$.   The graded module  associated is a
   module over $A[\eta_1,\dots,\eta_r]$.

   The radical of its annihilator ideal is independent of the choice of the generators $m_i$ so it defines
   a subvariety  $ch(M)$ of the cotangent bundle of  $T$ {\it the characteristic variety of $M$}.\smallskip
   
A basic fact  on $D-$modules  is that this variety is of dimension $\geq \dim M$ and, when it is of minimal dimension $=\dim M$, the module is said to be {\it holonomic}. It is again  a basic fact, see again \cite{Co}, that the coordinate ring of $\mathcal A_\Delta$  is a holonomic $D-$module and has a finite composition series. The task of the next section is to introduce and describe the irreducible modules which appear in this composition series as will finally be shown in  Proposition \ref{primadec}.  The special irreducible modules which appear can be thought of a the natural {\it quantization} of the cotangent bundles of the components of the arrangement. They have a very special description which allows us to put them together to produce explicit partial fraction decompositions  for the elements of the coordinate ring of the affine variety $\mathcal A_\Delta$ i.e. of the localized ring $\mathbb C[T][d^{-1}
],\ d:=\prod_{(a,\chi)\in\Delta}(1-a\chi).$
\subsection{Special Modules}
Take a sublattice $\Gamma\subset \Lambda$. Fix a maximal ideal $\mathfrak m\subset \mathbb
C[\Gamma]$ and denote by $I_{\mathfrak m}$ the extension of $\mathfrak m$ to $A$. The variety $W\subset T$ defined  by $I_{\mathfrak m}$ is a coset of the subgroup $T_\Gamma$ kernel of the characters in $\Gamma$ and will be called for short a {\it coset}.

$\Gamma$ is recovered from $W$ as the set of characters which are constant on $W$, clearly  $\mathfrak m$ is formed by the elements of  $\mathbb
C[\Gamma]$ vanishing on $W$.

Set $\Gamma^{\perp}=\{\phi\in\Lambda^*|\phi(\Gamma)=0\}$. Next let   $J_{\mathfrak m}$ be the left ideal in $D_T$ generated by $\mathfrak m$ and
by the derivations $D_{\phi}$ with $\phi\in \Gamma^\perp$.

{\bf Remark} The derivations $D_{\phi}$ with $\phi\in \Gamma^\perp$ are the invariant vector fields tangent to $W$.
\smallskip

We want to study  the structure of the $D_T$ module:
\begin{equation} N(W)=N(\Gamma, J_{\mathfrak m}):=D_T/J_{\mathfrak m}.
\end{equation}
We first treat the case in which $\Gamma$ is a
split direct summand, which is equivalent to saying that $T_\Gamma$ and $W$ are connected.

Under this assumption, we can then take a basis
 $\xi_1,\dots,\xi_r$ of the character group so that $\Gamma$ is generated by $\xi_1,\dots,\xi_k$  for some $k$.

 The derivations
 $\partial_i:=\xi_i{\partial\over \partial \xi_i}$  are a basis of the invariant vector fields. The $\partial_i,\ i=k+1,\dots,r$ are a basis of the space of  the derivations $D_{\phi}$ with $\phi\in \Gamma^\perp$.

 The ideal $\mathfrak m=(1-a_1\xi_1,\ldots ,1-a_k\xi_k)$ for  suitable $a_1,\ldots a_k\in \mathbb C^*$.

 \begin{lemma}\label{special}
 With the previous notations, and with $W$ a connected coset:

i) $N(W)$ is an irreducible $D_T$-module having as linear basis the classes of the elements $$t_{\underline h,\underline n}:=
\xi_{k+1}^{h_{k+1}}\dots \xi_r^{h_r}\prod_{ j=1,\dots,k}{\partial^{n_j}\over \partial \xi_j^{n_j}}$$$ h_s\in \mathbb Z,\ n_j\in \mathbb N$.

Or alternatively  the elements $$u_{\underline h,\underline n}:=\xi_{k+1}^{h_{k+1}}\dots \xi_r^{h_r}
\prod_{ j=1,\dots,k} \partial_j^{n_j} $$$ h_s\in \mathbb Z,\ n_j\in \mathbb N$.

 ii) $N(W)$ is a holonomic module with
characteristic variety the conormal bundle to $W$.
  \end{lemma}
  i) The  elements $\prod_{s=k+1}^r\xi_s^{h_s}\prod_{j=1}^k{\partial^{n_j}\over \partial \xi_j^{n_j}}\prod_{j=k+1}^r{\partial^{n_j}\over \partial \xi_j^{n_j}}\prod_{s=1}^k\xi_s^{h_s}$ with
  $h_s\in \mathbb Z$ and $\ n_j\in \mathbb N$ constitute  a basis for the ring $D_T$. This clearly implies that the classes of the $t_{\underline n,\underline h  }$'s  with $\ h_s\in \mathbb Z,\ n_j\in \mathbb N$, span $N(W)$.

   It remains  to show that the classes of the $t_{\underline n,\underline h  }$'s  with $\ h_s\in \mathbb Z,\ n_j\in \mathbb N$ are linearly independent
   and that $N(W)$ is irreducible.  Denote by $v$ the class of $1$ modulo $J_{\mathfrak m}$. Take an element  $$w=\sum_{(\underline
  h,\underline n)\in S}a_{\underline h,\underline n}t_{\underline h,\underline n}v$$ where $S$ is a subset of $\mathbb Z^{r-k}\times\mathbb N^k$
  and the $a_{\underline h,\underline n}$ are non zero complex numbers. Multiplying by a suitable monomial in $\xi_{k+1},\ldots ,\xi_r$, we can
  clearly assume that $S\subset \mathbb N^{r-k}\times\mathbb N^k=\mathbb N^r$.
  We want to show that the submodule $M$ generated by $w$ contains $v$. This will clearly implies both claims. We  make induction on the largest element $s$ of $S$ in the
  lexicographic ordering. If this element is zero there is nothing to prove since $w$ is a non zero multiple of $v$.

  Otherwise, write $s=(s_1,\ldots s_r)$ and let $1\leq i\leq r$ be the least element such that $s_i$ is non zero. If $i\leq k$, then we get
  $$(1-a_i\xi_i)w=\sum_{(\underline n,\underline h)\in S}a_{\underline h,\underline n}[(1-a_i\xi_i),t_{\underline h,\underline n}]v.$$
  Using the fact that $[(1-a_i\xi_i), {\partial\over \partial \xi_i}^{m}]=ma_i{\partial\over \partial \xi_i}^{m-1}$, we deduce that
  $$(1-a_i\xi)w=\sum_{(\underline n,\underline h)\in S'}a_{\underline h,\underline n}t_{\underline h,\underline n}v$$
  with the maximum element in $S'$ lexicographically smaller that the maximum element in $S$. Thus by induction we get that $v\
  $ lies in the $D_T$ submodule generated by $w$.

  If $i> k$, we apply ${\partial\over \partial \xi_i}$ and, reasoning in exactly the same way, we deduce that also in this case $M$ contains $v$.

For the other elements we use the fact that $\partial_i={\partial\over\partial \xi_i}\xi_i-1$ and see that we have a triangular base change.
\smallskip

ii) When we compute the characteristic variety we see that the graded module associated to $N(W)$ is still cyclic so it
equals $A[\eta_1,\dots,\eta_r]/J$ where $J$ is the ideal  generated by the elements  $1-a_j\xi_j, \ j=1,\dots,k,\ \eta_h  ,\ h=k+1,\dots,r$
which is the conormal bundle to $W$. \qed

Let us go back to a general $\Gamma\subset \Lambda$ and set $\overline \Gamma=\{\chi\in\Lambda|\exists\, n {\rm \  with\  }n\chi\in\Gamma\}$.
$\overline \Gamma$ is a split direct summand in $\Lambda$. Let $h=|\overline \Gamma/\Gamma|$.   The variety $W$ associated  to the ideal $I_{\mathfrak m}$ is now a union $W=\cup_{i=1}^hW_i$ of $h$ connected components, cosets of the
  connected torus $T_{\overline \Gamma}$.  These components are defined by the
maximal ideals   $\mathfrak  m_1,\ldots ,\mathfrak m_h\subset \mathbb
C[\overline\Gamma ]$ lying over $\mathfrak m$. We can clearly order the set $\{\phi_1,\ldots \phi_h\}$ of primitive
idempotents in $\mathbb C[\overline\Gamma ]/\mathfrak m\mathbb C[\overline\Gamma ]$ in such a way that the annihilator of $\phi_j$ in $\mathbb
C[\overline\Gamma ]$ is $\mathfrak  m_j$ for each $j=1,\ldots ,h$.

{\bf Remark} $\overline\Gamma$ is the set of characters in $\Lambda$ which are {\it locally constant} on $W$.
\begin{lemma}\label{induc} There is a canonical isomorphism of $D_T$ modules,
$$N(W)\simeq \oplus_{i=1}^hN(W_i).$$
\end{lemma}
\proof  Since clearly $J_{\mathfrak m}\subset J_{\mathfrak  m_i}$ for each $i=1,\ldots ,h$, we have a surjective homomorphism
$$f_i:N(W)\to N(W_i ).$$
Taking the direct sum of the $f_i$'s, we then get a homomorphism
$$f:N(W)\to \oplus_{i=1}^hN(W_i ).$$
Notice that since the $N(W_i )$'s are pairwise non isomorphic (having different characteristic varieties), we get that $f$ is surjective.

On the other hand, if we denote by $v$ the class of $1$ in $N(W)$ and set $v_i=\phi_iv$ (notice that this makes sense since
$\mathfrak m$ annihilates $v$), we have that $J_{\mathfrak m_i}$ annihilates $v_i$. We thus get  homomorphisms
$$g_i:N(W_i )\to N(W)$$
$i=1,\ldots h$. Taking their  sum we get
$$g:\oplus_{i=1}^h N(W_i )\to N(W)$$
Since $v=\sum_{i=1}^hv_i$, we immediately deduce that $g$  is surjective.

 We  leave to the reader the simple verification  that $f$ and $g$ are one inverse to the other.
\qed

We are now going to give explicit realizations of the module $N(W)$.

Before doing this, we need a general statement which is well known but whose proof we include for completeness. Let $R$ be a commutative ring.
$\underline z=(z_1,\ldots ,z_k)$ a sequence of elements  in $R$ such that any permutation of them is a regular sequence. For any subset
$S\subset \{1,\ldots ,k\}$, set $z_S=\prod_{j\notin S } z_j$. Consider the \v{C}ech  complex $C(R,\underline z)$
$$0\to R\to\oplus_{i=1}^kR[z_i^{-1}]\to\cdots\to \oplus_{S,|S|=t}R[z_S^{-1}] \to \cdots \to R[(z_1\cdots z_k)^{-1}]\to 0$$

\begin{lemma}\label{local} 1) $H^i(C(R,\underline z))=0$ for $i=0,\ldots k-1$.

2) The $R$ submodule of $H^k(C(R,\underline z) ) $, generated by the class of $(z_1\cdots z_k)^{-1}$ is isomorphic to $R/(z_1,\ldots ,z_k)$. In particular
$H^k(C(R,\underline z) )\neq 0$.
\end{lemma}
\proof We proceed by induction on $k$. If $k=1$ our complex is
$$0\to R\to R[z^{-1}]\to 0$$
the map $R\to  R[z^{-1}]$ is injective since $z$ is not a zero divisor. It is not  surjective since $z$ is not invertible. Clearly the
annihilator of the class of $1/z$ in $R[z^{-1}]/R$ consists of those elements $r\in R$ such that $r/z\in R$. This is  the  ideal $(z)$.

We now assume  our statement  for $k-1$. Consider the ring $R'=R/(z_k)$ and in it the  sequence $\underline z'=(z'_1,\ldots ,z'_{k-1})$ where
$z'_j$ is the image of $z_j$ modulo $z_k$. Our  assumptions are satisfied by $R'$ and the sequence $\underline z'$. Also we have an exact
sequence of complexes
$$0\to C(R,\underline z)\stackrel{\cdot z_k}{\to} C(R,\underline z)\to C(R',\underline z')\to 0.$$
Applying the inductive assumption and the long cohomology sequence we immediately deduce that $H^i(C(R,\underline z))=0$ for $0\leq i\leq k-2$
and multiplication by $z_k$ gives an injective map
$$0\to H^{k-1}(C(R,\underline z))\stackrel{\cdot z_k}{\to}H^{k-1}(C(R,\underline z)).$$
We claim that $z_k$ acts locally nilpotently on $ H^{k-1}(C(R,\underline z))$ so that necessarily
 $H^{k-1}(C(R,\underline z))=0$ Indeed let $(r_1,\ldots ,r_k)\in \oplus _{i=1}^kR[z_{\{i\}}]$ be a cycle.
 This means that $r_1+r_2+\cdots +r_k=0$. If we multiply by a big enough power of $z_k$, say $z_k^M$,
 we get that $z_k^Mr_i\in  R[z_{\{i,k\}}]$ for each $i=1,\ldots ,k-1$. Then $z_k^M(r_1,\ldots ,r_k)=
( z_k^Mr_1,\ldots ,z_k^Mr_k)$ is cohomologous to $(0,\ldots ,0, z_k^M(r_1+r_2+\cdots +r_k))= (0,\ldots ,0)$ as desired.

It remains to show 2). For this again using the long cohomology sequence and 1) we get an exact sequence
$$0\to H^{k-1}(C(R',\underline z')\stackrel{\delta}{\to}
H^{k}(C(R,\underline z))\stackrel{\cdot z_k}{\to} H^{k}(C(R,\underline z))\to 0.$$ It is easy to see that $\delta((z'_1\cdots
z'_{k-1})^{-1})=(z_1\cdots z_{k})^{-1}$. This and our inductive assumption immediately imply our claim (notice that $R/(z_1,\ldots
,z_k)=R'/(z'_1,\ldots ,z'_{k-1})$).\qed
\begin{remark} It is well known that $H^i(C(R,\underline z))$ equals the $i$-th local cohomology group $H^i_{(z_1,\ldots z_{k})}(R) $
(see for example \cite {Ha} where our Lemma \ref{local} is essentially contained).\end{remark}

We are now ready to give explicit realizations of the module $N(W)=N(\Gamma; J_{\mathfrak m})$. To achieve this let us choose a basis
$\Psi:=\{\psi_1,\dots,\psi_k\}$ in $\Lambda$ of $\Gamma$ so that the ideal $\mathfrak m=(1-a_1\psi_1,\ldots ,1-a_k\psi_k)$ for suitable non zero constants
$a_i\in \mathbb C^*$ defines $W$.  Set $d=\prod_{i=1}^k(1-a_i\psi_i)$. and for any $i=1,\ldots ,k$, $d_i=d/(1-a_i\psi_i)$. Consider the ring $A[d^{-1}]$ as
$D_T$ module and remark that $A[d_i^{-1}]$ is a submodule for each $i$. We have
\begin{proposition}\label{decomp0} The $D_T$ module $A[d^{-1}]/\sum_{i=1}^kA[d_i^{-1}]$ is isomorphic to $N(W)$.
\end{proposition}
\proof As before let  $\mathfrak m_1,\ldots ,\mathfrak m_h\subset \mathbb C[\overline\Gamma ]$ be set set of maximal ideals lying over $\mathfrak m$,
 the set $\{\phi_1,\ldots \phi_h\}$  the set  of primitive idempotents in $\mathbb C[\overline\Gamma ]/m\mathbb C[\overline\Gamma ]$
  ordered in such a way that the annihilator of $\phi_j$ in $\mathbb C[\overline\Gamma ]$ is $\mathfrak m_j$ for each $j=1,\ldots ,h$.

Denote by $w$ the class of $d^{-1}$ in $A[d^{-1}]/\sum_{i=1}^kA[d_i^{-1}]$. It is immediate to see that $J_{\mathfrak m}$ annihilates $w$ so that we
get a map of $D_T$ modules
$$f:N(\Gamma; J_{\mathfrak m})\to A[d^{-1}]/\sum_{i=1}^kA[d_i^{-1}]$$
with $f(v)=w$. To see that this map is injective it suffices to show that its restriction to each irreducible component $N(\overline\Gamma;
J_{\mathfrak m_i})$ is non zero. $N(\overline\Gamma; J_{\mathfrak m_i})$ is generated by $v_i=\phi_iv$.
 We have $f(v_i)=\phi f(v)=\phi w$ and this is non zero by lemma \ref {local}.

Thus we only have to show that $f$ is surjective. $A[d^{-1}]/\sum_{i=1}^kA[d_i^{-1}]$ is spanned by the classes of the elements
$${a\over (1-a_1\psi_1)^{n_1}\cdots (1-a_k\psi_k)^{n_k}}$$
with $a\in A$ and $n_i>o$ for each $i=1,\ldots ,k$. In $\Lambda^*\otimes \mathbb Q$ choose elements $\rho_1,\ldots \rho_k$ such  that $\langle
\rho_j,\psi_i\rangle=\delta_{i,j}$ for each $i,j=1,\ldots ,k$. Take the derivations $D_i:=(a_i\psi_i)^{-1}D_{\rho_i}$. It is then immediate to
see that
$${aD_1^{n_1}\cdots D_k^{n_k}\over n_1!\cdots n_k!}\Big ({1\over (1-a_1\psi_1)\cdots (1-a_k\psi_k)}\Big )={a\over (1-a_1\psi_1)^{n_1}\cdots (1-a_k\psi_k)^{n_k}}$$
and everything follows.\qed

\section{Toric Arrangements}
\subsection{Toric arrangements}
Let $T=hom_{\mathbb Z}(\Lambda,\mathbb C^*)$ be an algebraic torus over $\mathbb C$, of dimension $r$. $\Lambda$ its character group,   a rank $r$ lattice,
$A:=\mathbb C[T]=\mathbb C[\Lambda]$ its coordinate ring.

We want to study in this setting a suitable analogue of the theory of hyperplane
arrangements. We develop the theory in a level of generality higher than the one needed for the applications to counting formulas. \smallskip

 We consider, as analogue of a hyperplane, a translate of the kernel of a character $\chi$, i.e. the hypersurface
$H_{a,\chi}$ of equation $1-a\chi=0$ for a suitable non zero $a\in \mathbb C^*$. \smallskip

Notice that, if $\chi$ is not primitive, this hypersurface is not connected, but it is the union of finitely many cosets of the codimension 1 torus, the connected component of the identity of  the kernel of $\chi$.
If $\Psi:=\{\chi_1,\dots,\chi_i\}$ are linearly independent characters and $\underline a=(a_1,\ldots ,a_i)$ is a sequence of non zero complex
numbers,  the intersection $K_{\underline a,\Psi}$ of the $H_{a_s,\chi_s}$ is in general not connected even if the characters are primitive.

 In fact let us denote by  $\langle\Psi\rangle:=\langle \chi_1,\dots,\chi_i \rangle$ the subgroup generated by the $\chi_j,j=1,\dots,i$,  and  by
$\overline{\langle\Psi\rangle}$ the subgroup of $\Lambda$ formed by the elements $\chi$ such that $\chi^m\in \langle\Psi\rangle$ for some $m>0$. The following is immediate:

\begin{proposition}  $K_{\underline a,\Psi}$ is a coset of the subgroup kernel of the characters in $\langle\Psi\rangle$.
It has 
$|\overline{\langle\Psi\rangle}/\langle\Psi\rangle|$   connected components  each of which is a coset of the subgroup kernel of the characters in $\overline{\langle\Psi\rangle}$.
 \end{proposition}\medskip

Given a finite subset  $\Delta\subset \mathbb C^*\times\Lambda$, we shall study the components of the  {\bf arrangement $\mathcal R_\Delta$} generated by the
hypersurfaces $H_{a,\chi}$ as $(a,\chi)\in\Delta$. 

As explained in the introduction, by this we shall mean {\it the set of all connected components} of all intersections of these
hypersurfaces. \smallskip

The complement of the union of these hypersurfaces will be denoted by $\mathcal A_\Delta$, it is an affine variety of coordinate ring $A[d^{-1}
],\ d:=\prod_{(a,\chi)\in\Delta}(1-a\chi).$  Our main concern will be to compute the cohomology of $\mathcal A_\Delta$. 
\subsection{The $D_T$-module $R$}  With the notations of the previous section,
if $W$ is one component of the arrangement $\mathcal R_\Delta$, let $\Delta(W)$ be the set of elements $(a,\chi)\in\Delta$ with $1-a\chi=0$ on $W$ and $\Sigma(W)$ the set of characters $\chi\in \Lambda$ which, as functions on $W$, are constant. It is easily seen that $\Sigma(W)$ is a split direct summand on $\Lambda$ and that $W$ is a coset of the subtorus   $T(W)$ kernel of the characters of $\Sigma(W)$.\smallskip

The subtori   $T(W)$ will be called the {\it tori associated to the arrangement}.\smallskip

The complement of the union of the hypersurfaces $H_{a,\chi}$ will be denoted by $\mathcal A_\Delta$, it is an affine variety of coordinate ring
$R=A[d^{-1} ],$ with $ d:=\prod_{(a,\chi)\in\Delta}(1-a\chi).$  We begin by studying the ring $R$ as a $D_T$ module.
\smallskip

Let  $\pi:\Delta \to
\Lambda$
 be  the projection on the second
factor and set $\overline \Delta:=\pi(\Delta)$.

We shall say that a subset $S=  \{(a_1,\psi_1),\dots, (a_r,\psi_r)\}\subset \Delta$ is linearly independent if  $ \{ \psi_1 ,\dots,  \psi_r \}\subset \Lambda$ are linearly independent. We set for each $k=0,\ldots ,r$, $$\mathcal J_k=\{S\subset
\Delta|\ S  {\rm \ is\  linearly\  independent\  and\ } |S|\leq k\}.$$ Of course $\mathcal J_{r+s}=\mathcal J_r:=\mathcal J$ for each $s\geq
0$ ($J_0$ is by definition the empty set). Given $S\in \mathcal J$ we set $d_S:=\prod_{(a,\chi)\in S}(1-a\chi)$ ($d_{\emptyset}=1$), and define
for each $k=0,\ldots ,r$,
$$R_k=\sum_{S\in\mathcal J_k}A[d_S^{-1}].$$
Clearly $R_k\subset R_{k+1}$ for each $k$ and $R_k$ is a sub $D_T$ module of $R$.

The next Proposition  is a special case of Proposition 4.6 of \cite{S1}.
 \begin{proposition}\label{riduz} Let  $\Psi\subset
 \mathbb C^*\times\Lambda$ be  finite. Let
 $\Xi\subset \Lambda$ denote the projection of
 $\Psi$ to the second factor and set $\Gamma$
 equal to  the lattice  generated by  $\Xi$.
 Then $$\prod_{(a,\psi)\in\Psi}{1\over (1-a
 \psi)}$$
can be written as a linear  combination with coefficients in  the ring of $\mathbb C[\Gamma]$ of elements of the form $${1\over (1-
a_1\psi_{1})^{h_1}\cdots (1- a_r\psi_{r})^{h_r}}$$  where    $(a_i,\psi_i)\in \Psi$ for each $i=1,\dots ,r$ and $\{\psi_1,\dots,\psi_r\}$ are
linearly independent.
\end{proposition}
\proof
 By  a simple induction we  may assume that $\Psi=\{(a_0,\psi_0),\ldots ,(a_n,\psi_n)\}$ with $\psi_0,\ldots \psi_n$ linearly dependent and
$\psi_1,\ldots \psi_n$  linearly independent.   We then choose  representatives $\phi_i\in \mathbb C[\Gamma]$ of the primitive idempotents  $e_i$  of the ring
$\mathbb C[\Gamma]/(1-a_1\psi_1,\dots,1-a_n\psi_n)$. Notice that  in this ring we have  $\psi_0e_i=\beta_ie_i$ with $\beta_i\in \mathbb C$.  By
the definition of the $\phi_i$'s,   we have $1=\sum_i\phi_i+\sum_{j=1}^nb_i(1-a_j\psi_j)$, for some $b_i\in \mathbb C[\Gamma]$. So $${1\over (1-a_0 \psi_0) \cdots (1- a_n\psi_n)
}={\sum_i\phi_i+\sum_{j=1}^nb_i(1-a_j\psi_j)\over (1- a_0\psi_0) \cdots (1- a_n\psi_n) }=$$ $$\sum_i {\phi_i \over (1-a_0 \psi_0) \cdots (1-
a_n\psi_n)  }+\sum_{j=1}^n{b_i(1-a_j\psi_j)\over (1-a_0 \psi_0) \cdots (1- a_n\psi_n) }$$ Using induction on the cardinality of $\Psi$ we then
have to  analyze only the terms $$ { \phi_i \over (1-a_0 \psi_0) \cdots (1- a_n\psi_n)}.$$ We have  that $(1-a_0\psi_0)\phi_i=\sum_{j=1}^nc_j(
1-a_j\psi_j)+\gamma_i\phi_i$  with $\gamma_i=1-a_0\beta_i.$ We separate two cases. If  $\gamma_i=0$, we have
$(1-a_0\psi_0)\phi_i=\sum_{j=1}^nc_j( 1-a_j\psi_j)$ and  substituting $$ { \phi_i \over (1-a_0 \psi_0) \cdots (1- a_n\psi_n) }=   {
\sum_{j=1}^nc_j( 1-a_j\psi_j)\over (1- a_0\psi_0)^2(1-a_1\psi_1) \cdots (1- a_n\psi_n) }$$  we obtain a sum of terms in which in the denominator
some $1-a_j\psi_j$ has disappeared.

If $\gamma_i\neq 0$, we get $${\phi_i\over (1-a_0 \psi_0) \cdots (1- a_n\psi_n) }=\gamma_i^{-1} {(1-a_0\psi_0)\phi_i-\sum_{j=1}^nc_j(
1-a_j\psi_j)\over (1- a_0\psi_0) \cdots (1- a_n\psi_n) },$$ and again everything follows by induction.\qed

\subsection{Partial fractions for $R$}
From Proposition \ref{riduz} we
deduce
\begin{corollary} $R=R_r$.
\end{corollary}
\proof The proof is clear.\qed

Let us now take $S\subset \mathcal J_k$. Let $\Gamma \subset \Lambda$ be the sublattice generated by the characters $\pi(S)$, $\mathfrak m$ the
maximal ideal in $\mathbb C[\Gamma]$ generated by the elements $1-a\chi$ with $(a,\chi)\in S$. As in the previous section introduce the lattice
$\overline \Gamma=\{\chi\in\Lambda|\exists\, n {\rm \  with\  } \chi^n\in\Gamma\}$. Set $h=|\overline \Gamma/\Gamma|$. Let $\mathfrak m_1,\ldots , \frak
m_h\subset \mathbb C[\overline\Gamma ]$ be the set of maximal ideals lying over $\mathfrak m$, and order the set $\{\phi_1,\ldots \phi_h\}$ of
primitive idempotents in $\mathbb C[\overline\Gamma ]/m\mathbb C[\overline\Gamma ]$ in such a way that the annihilator of $\phi_j$ in $\mathbb
C[\overline\Gamma ]$ is $\mathfrak m_j$ for each $j=1,\ldots ,h$.

 From Proposition \ref{decomp0} we deduce  a canonical isomorphism  as $D_T$-modules:
\begin{equation}
\label{iso} A[d_S^{-1}]/(R_{k-1}\cap A[d_S^{-1}])\simeq N(\Gamma; J_{\mathfrak m})\simeq \oplus_{i=1}^hN(\overline\Gamma; J_{\frak
m_i})\end{equation} mapping the class of $1$ in  $N(\overline\Gamma; J_{\mathfrak m_i})$ to $\phi_i[1/d_S]$, with $[1/d_S]$ the class of the vector $1/d_S$.

\medskip

Let $W$ be a component of the arrangement of codimension $k$. We shall say that  $S\subset \mathcal
J_k$ and $W$ are associated $W$ is a component of the variety defined by the vanishing of the elements $1-a\chi$ with $(a,\chi)\in S$.  In this case denoting by $\mathfrak n_S$ the ideal generated by the elements $1-a\chi $ with $(a,\chi)\in S$ we have that $\mathbb C[\Lambda]/\mathfrak n_S$ has a primitive idempotent $\phi_{W,S}$ which, as function on $Z_S$  is the characteristic function of the component $W$.

We  associate to such a $W$  a subspace $V_{W}\subset R_k/R_{k-1}$ as follows. 
\begin{definition} Define $V_{W}$ as
the subspace spanned by the classes $v_{S,W}:=\phi_{W,S }[1/d_S]$ as $S$ runs among the subsets of  $\mathcal  J_k$
 associated to
$W$.\end{definition}
 Notice that the $\mathbb C[\Lambda]$ module generated by $d_S^{-1}$ in $R_k/R_{k-1}$ factors as an $\mathbb C[\Lambda]/\mathfrak n_S$ module.
\begin{proposition}\label{primadec} For each $k=0,\ldots r$, $R_k/R_{k-1}$ is a semisimple $D_T$ module and we have a canonical isomorphism
\begin{equation}\label{isom} R_k/R_{k-1}\simeq \oplus_{W} N(W)\otimes V_{W}\end{equation}
$W$ runs over the components of the arrangement of codimension $k$.\end{proposition}
\proof  $R_k/R_{k-1}$ is a quotient of $\oplus_{S\in\mathcal J_k}A[d_S^{-1}]/(R_{k-1}\cap A[d_S^{-1}])$. From the fact that each
$A[d_S^{-1}]/(R_{k-1}\cap A[d_S^{-1}])$ is semisimple the semisimplicity of $R_k/R_{k-1}$ follows. The rest of the Proposition follows from
description of the canonical isomorphism (\ref {isom}) recalled above.\qed

 \bigskip

 \subsection{No broken circuits}  In order to finish our analysis of $R$ as a $D_T$ module we have to understand the dimension of $V_{W}$ for each component  $W$ or better yet, to give a basis for these spaces.  From now on   we shall choose a component $W$.

 Fix a total order on the set of characters $\overline\Delta$.  Set
 $$\Delta_{W}=\{(a,\chi)\in\Delta|\ \chi\in\Sigma {\rm \ and \ }1-a\chi=0,\ \text{on}\ W\}.$$
 The projection $\pi$ (cf. 3.1) of $\Delta$ to $\overline\Delta$ restricted to $\Delta_{W}$ is injective  thus 
 the total ordering on $\overline\Delta$ induces a total ordering on $\Delta_{W}$.

       \begin{definition} A subset $S=\{(a_1,\chi_1)<\cdots <(a_k,\chi_k)\}$ of $\Delta_{W}$ is called a non broken circuit on $W$ (relative to our chosen ordering) if
\begin{enumerate} \item $S$ is associated to $W$
 \item There is no    $(a,\chi)\in \Delta_{W}$ and no $e=1,\ldots ,k$ with the property that
   $(a,\chi)<(a_e,\chi_e)$ and $\chi,\chi_{e}, \ldots ,\chi_{k}$ are linearly dependent characters.
\end{enumerate}
  \end{definition}

  \begin{theorem}\label{nobroken} The classes $v_{S,W}$, as $S$ runs among the non broken circuits in $\Delta_{W}$ form a basis of
 $V_{W}$.
\end{theorem}
\proof We first show that the classes of $v_{S,W} $ with $S$ a non broken circuit, span $V_{W}$. For this let us order
the $S$'s associated to $W$ lexicographically. It is clear that the minimum $S$ in this ordering is non broken.

Pick a general $S=\{(a_1,\chi_1)<\cdots <(a_k,\chi_k)\}$ associated to $W$ and assume there is     $(a,\chi)\in
\Delta_{W}$ with
  $\chi,\chi_{e}, \ldots ,\chi_{k}$  linearly dependent characters,
  and $(a,\chi)<(a_e,\chi_e)$. We can also assume, by induction, that if we remove one of the $\chi_i,  e\leq i\leq k$ the remaining elements  $\chi,\chi_{e}, \ldots ,\check \chi_i,\ldots ,\chi_{k}$ are  linearly independent.
  
    Clearly   for each $j=0,\ldots ,k-e$,
   the set $S_j=(S-\{(a_{e+j},\chi_{e+j})\})\cup \{(a,\chi)\}$ is associated to $W$.

Let $\Gamma\supset \Gamma'$ be the sublattices generated by $\chi,\chi_{e}, \ldots ,\chi_{k}$ and $  \chi_{e}, \ldots ,\chi_{k}$ respectively. Let $U$ be the torus of coordinate ring $\mathbb C[\Gamma]$. We have an associated map $\rho:T\to U$ under which $W$ maps to the point $P_0$ of equations $1-a\chi=1-a_e\chi_e=\ldots=1-a_k\chi_k=0$.

If we consider only the equations $1-a_e\chi_e=\ldots=1-a_k\chi_k=0$, they define finitely many points $P_0,P_1,\dots,P_t$ in $U$, we can
thus choose an element $\psi\in\mathbb C[\Gamma]$  which takes value 1 on $P_0$ and 0 in the remaining $P_i$.  In particular $\psi=1$ on $W$ and
 $$  \phi_{S,W}=\psi\phi_{S,W},\quad v_{S,W} =\psi v_{S,W}.$$
We have $\psi(1-a\chi)$ vanishes on all the points $P_0,P_1,\dots,P_t$ hence  $\psi(1-a\chi)=c_e  (1-a_e\chi_e)+ \ldots +c_k( 1-a_k\chi_k)$, with $c_i\in\mathbb C[\Gamma]$, hence:
$${\psi\over \prod_{i=e}^k(1-a_i\chi_i)}={\psi(1-a\chi)\over (1-a\chi)
\prod_{i=e}^k(1-a_i\chi_i)}=\sum_{j=e}^{k }{c_j\over (1-a\chi)\prod_{i=1, i\neq j }^k(1-a_i\chi_i)}.$$
Choose a representative $f_{S,W}$  of $\phi_{S,W}$, we have from the previous identity:
\begin{equation}\label{linc}v_{S,W} =\sum_{j=e}^{k }[{f_{S,W}c_j\over \prod_{i=1}^{e-1}(1-a_i\chi_i)(1-a\chi)\prod_{i=e, i\neq j }^k(1-a_i\chi_i)} ].\end{equation}
We have to understand the elements $$[{f_{S,W}c_j\over \prod_{i=1}^{e-1}(1-a_i\chi_i)(1-a\chi)\prod_{i=e, i\neq j }^k(1-a_i\chi_i)}]\in R_k/R_{k-1}.$$
In $U$, for each $j$ between $e$ and $k$ the equations $1-a\chi=0, 1-a_i\chi_i=0,\ \forall e\leq i\leq k, \ i\neq j$  define a finite set of points $U_j$.  Let $Q:=\cup U_j=\{P_0, A_1,\dots,A_f\}$ again a finite set of points including $P_0$.  We can assume that we have chosen   $f_{S,W}$ with the further property of vanishing on all the subvarieties $\rho^{-1}(A_j)$. Consider next the subvariety $Z_j$ of $T$ of equations $1-a\chi=0, 1-a_i\chi_i=0,\ \forall 1\leq i\leq k, \ i\neq j$, its connected components $W_h$ are elements of the arrangement. $Z_j$ is the preimage under the map $\rho:T\to U$ of the finite number of points given by the same equations in $U$.  We claim that on each such $W_h$ the function $f_{S,W}c_j$ is constant.  This   implies  that, modulo the ideal generated by the elements $1-a\chi=0, 1-a_i\chi_i=0,\ \forall 1\leq i\leq k, \ i\neq j$ the element $f_{S,W}c_j$ is a linear combination of primitive idempotents, hence  by (\ref{linc}) $v_{S,W}$  is a linear combination of elements $v_{W_h,S_i}$ with $S_i$   lower than $S$ in the lexicographic order.

In fact this is clear for $c_j$ since this comes from a function on $U$, as for  $f_{S,W}$ we know, by hypothesis that it vanishes on all the components $W_h$ which map  under $\rho$ to a point $A_p,\ p=1,\dots,f$ thus it only remains to analyze the components lying over $P_0$. On these components also $1-a_j\chi_j$ vanishes and so they are part of the components of the subvariety  of equations $1-a_i\chi_i=0, \ 1\leq i\leq k$. By construction, on this subvariety,  $f_{S,W}$ is the characteristic function of $W$, hence the claim. From this, everything follows by induction.
\medskip

We now need to show that the vectors $v_{S,W}$ with $S$ a non broken circuit are linearly independent. Take a basis $\xi_1,\ldots ,
\xi_k$ of $\Sigma$ and complete it to a basis $\xi_1,\ldots , \xi_r$ of $\Lambda$. Then $I_{\mathfrak n}=(1-b_1\xi_1,\ldots ,1-\xi_k)\subset A$ for
suitable $b_j\in\mathbb C^*$. Consider the completion $B$ of the local ring $A_{I_{\mathfrak n}}$ at  its maximal ideal. We can identify $B$ with
the power series ring $K[[x_1,\ldots ,x_k]]$ with $x_j=\log(b_j\xi_j)$ and $K=\mathbb C(\xi_{k+1},\ldots ,\xi_{r})$, the quotient field of
$A/I_{\mathfrak n}$. From now on we shall identify $A$ with its image under its canonical inclusion in $B$. Thus $b_j\xi_j$ equals $\exp
(x_j)=\sum_{s}{x_j^s/ s!}.$ Similarly  for $(a,\chi)\in \Delta_{W}$, $a\chi$ equals $\exp (z_{\chi})$ where if
$\chi=\xi_i^{m_1}\cdots \xi_k^{m_k}$, $z_{\chi}=m_1x_1+\cdots +m_kx_k$, in particular this implies that, in $B$, we have  $1-a\chi=z_{\chi}f_{\chi}$ with
$f_{\chi}$ invertible in $B$ and congruent to 1 modulo the maximal ideal. On the other hand, if $(a,\chi)\notin \Delta_{W}$ then
$1-a\chi$ is invertible in $B$. It follows that
$$\hat R=B[\prod _{(a,\chi)\in \Delta}{1\over (1-a\chi)}]=B[\prod _{(a,\chi)\in \Delta_{W}}{1\over z_{\chi}}]$$
Let us define a filtration of  $\hat R$ in a completely analogous way to the filtration of $R$. Notice that $R_s\subset \hat R_s$ for each $s$,
so that we get a map $j:R_k/R_{k-1}\to \hat R_k/\hat R_{k-1}$. The above considerations clearly imply that for any  $S=\{(a_1,\chi_1)<\cdots
<(a_k,\chi_k)\}$ associated to $W$ we have that $j(v_{S,\mathfrak n})$ equals the class  $w_S$ of $1/(z_{\chi_1}\cdots z_{\chi_k})$.
Using this, our claim follows immediately from the linear independence the $w_S$'s  as $S$ runs among the non broken circuits which is well
known (see \cite{te}).\qed
\subsection{Partial fractions}

    Keeping the notations of the previous section consider a component $W$  of the arrangement of codimension $k$. Its associated character group   $\Sigma(W):=\Sigma$ is  a split direct summand of $\Lambda$ of rank $k$ and $W$ is defined by the ideal generated by
     a maximal ideal  $\mathfrak n$   in $\mathbb C[\Sigma]$. By assumption there exists at least a subset $S\subset \Delta$ associated to
    $W$. This means that the irreducible subvariety $W$ is a connected component of the intersection of the
    $H_{(a,\chi)}$ with $(a,\chi)\in S$.\smallskip

    Take a basis $\xi_1,\ldots ,
\xi_k$ of $\Sigma(W)$ and complete it to a basis $\xi_1,\ldots , \xi_r$ of $\Lambda$. For each $i=1,\ldots r$, set $\partial_i=\partial/\partial
\xi_i$. This allows us to introduce a  commutative subring in $D_T$ namely $D_{W}:=\mathbb C[\xi_{k+1}^{\pm 1},\ldots ,\xi_{r}^{\pm
1}][\partial_{1},\dots ,\partial_k]$.

In $R_k/R_{k-1}$ consider the  component $N(W)\otimes V_{W}$. Let us remark that Lemma \ref{special} implies
that as a $D_{W}$-module $N(W)$ is  free of rank 1, generated by the class of 1.

Also remark that by Theorem \ref{nobroken}, we have a basis of $V_{W}$ given by the vectors $ v_{S,W} $
 as $S$ runs through the non broken circuits in $\Delta_{W}$.
 We choose a representative $f_{S,W}/\prod_{(a,\chi)\in S}(1-a\chi)$
 of  $v_{S,W} $ where  $f_{S,W}$ is a representative of the primitive idempotent $ \phi_{S,W} $.  In this way we   lift
  $V_{W}$ to a
subspace $\tilde V_{W}$ lying in the ring obtained from $\mathbb C[\Sigma]$ inverting the product of the
elements $1-a\chi$ with $ (a,\chi)\in\Delta_{W}$.      Let $M(W)$ be the $D_W$-module generated by  $\tilde V_{W}$,  $M(W)$ is a free  $D_W$-module with basis any linear basis of $\tilde V_{W}$, in particular the elements $f_{S,W}/\prod_{(a,\chi)\in S}(1-a\chi)$. Under the quotient map  $R_k\to R_k/R_{k-1}$,   $M(W)$ maps isomorphically   to $N(W)\otimes V_{W}$. In particular we get,

 \begin{corollary}\label{direct}   $R =\oplus_{W}  M(W)$.

       If the codimension of $W$ is $k$, $M(W)\subset R_k$ and $\pi_k$ maps $M(W)$ isomorphically onto
         $N(W)\otimes V_{W}$.
         \end{corollary}\label{dec}

 We can make explicit the previous decomposition in terms of an explicit linear basis and a suitable choice of representatives of the primitive idempotents.

 For every non broken linearly independent subset $\Gamma=\{\psi_1,\dots,\psi_k\}\subset\Delta$  let us denote by $\Lambda(\Gamma),\overline \Lambda(\Gamma)$ the lattice that $\Gamma$ generates and the set of elements of $\Lambda$ which are torsion modulo $\Lambda(\Gamma)$.   
The idempotents $\phi_{ \Gamma,W}$ of the ring  $\mathbb C[\Lambda]/\langle 1-a_1 \psi_1,\dots,1-a_k\psi_k\rangle$ can be also described in terms of characters of the finite group $\overline \Lambda(\Gamma)/\Lambda(\Gamma).$ 
First if all the $a_i=1$,
if $\Phi(\Gamma)$ denotes the character group, given $\lambda\in \Phi(\Gamma)$ a representative $\phi_\lambda$ of the corresponding idempotent  is constructed as follows, let $h_\Gamma:= |\overline \Lambda(\Gamma)/\Lambda(\Gamma)|$,  take representatives
 $\xi_1,\dots,\xi_h\in \overline \Lambda(\Gamma) $ for the cosets of $\overline \Lambda(\Gamma)/\Lambda(\Gamma)$ then:
 \begin{equation}\label{kar} \phi_\lambda={1\over h_\Gamma}\sum_{i=1}^{h_\Gamma} \lambda(\xi_i^{-1})\xi_i .\end{equation}
 We can
further normalize the $x_i$  as follows. We think of $\Lambda$ as being contained n the rational vector space with basis the elements $\psi_i$  then we can uniquely choose the $x_i$ of the form $\prod_j\psi_j^{a_{ij}}$ with $0\leq a_{ij}<1$  rational numbers. In this way $\phi_\lambda$ is a well defined element of $R$.

  When the $a_i\neq 1$ we must make a change of variables in the ring $\mathbb C[\Lambda]$ so that $a_i\psi_i$ becomes one of the group elements (apply an automorphism of translation). We then have
  \begin{theorem}\label{gcase}   Every element   $f\in R$ can be uniquely expressed in a {\it decomposition in partial fractions}:
      \begin{equation} \label{gdec} f= \sum_\Gamma \sum_{\lambda\in\Phi(\Gamma)}\sum_{n_1,\dots,n_k} { a_{n_1,\dots,n_k}^{x_i,\Gamma} x_i\over (1-a_1\psi_1)^{n_1}\dots(1-a_k\psi_k)^{n_k}}  \end{equation}
As $x_i$ runs over a set of representatives of the cosets of $\overline \Lambda(\Gamma)/\Lambda(\Gamma) $ or as sum:
      \begin{equation} \label{gdeco} f= \sum_\Gamma \sum_{\lambda\in\Phi(\Gamma)}\sum_{n_1,\dots,n_k} { a_{n_1,\dots,n_k}^{\lambda,\Gamma} \phi_\lambda\over (1-a_1\psi_1)^{n_1}\dots(1-a_k\psi_k)^{n_k}}  \end{equation}
  as $\Gamma$ runs over all the no broken circuit sets and $a_{n_1,\dots,n_k}^{\lambda,\Gamma}\in  \mathbb C[\xi_{k+1}^{\pm 1},\ldots ,\xi_{r}^{\pm
1}]$ in the corresponding ring.   \end{theorem}
  \proof
  The pairs $\Gamma,\lambda\in \Phi(\Gamma)$ are in $1-1$  correspondence with the components $W$ of the arrangement and we have that the elements ${\xi_{k+1}^{h_1}\dots \xi_r^{h_r}\phi_\lambda\over (1-a_1\psi_1)^{n_1}\dots(1-a_k\psi_k)^{n_k}},$ $\ n_i\in\mathbb  N,\ h_i\in\mathbb  Z$ form a basis of the corresponding vector space $M(W)$.  Summing these spaces over the components associated to $\Gamma$ we have a free module over the ring $D_{W}:=\mathbb C[\xi_{k+1}^{\pm 1},\ldots ,\xi_{r}^{\pm
1}][\partial_{1},\dots ,\partial_k]$  with basis the representatives $x_i$. Then the formula is clear.
 
  \qed
\section{Cohomology}
    In this section we shall determine the additive structure of the cohomology of the open set $\mathcal A_{\Delta}\subset T$ which is the
    complement of the union of the hypersurfaces of $H_{(a,\chi)}$ of equation $1-a\chi=0$ for $(a,\chi)\in \Delta$.

Let   $\Xi=\bigwedge(\omega_1,\dots,\omega_r)$ denote the exterior algebra of invariant differential forms, $\omega_i:=d\log \xi_i$ for our
basis of characters.

Due to the Theorem of Grothendieck \cite{gr}, the De Rham cohomology of $\mathcal A_{\Delta}$ can be computed via its algebraic
 De Rham complex    $\Omega:=R\otimes \Xi$.

It is convenient to write the differential in an invariant way.
$$ df=\sum_i{\partial f\over \partial \xi_i}d\xi_i=\sum_i\partial_i(f)\omega_i.$$
Using the decomposition of Corollary \ref{direct}, we also have:
$$\Omega =\oplus_{W} M(W) \otimes\Xi.$$
Next we claim that $ M(W) \otimes\Xi$ is a subcomplex.

For this it is enough to see that each  $D_W f_{S,W}/\prod_{(a,\chi)\in S}(1-a\chi) \Xi$ is a complex.
Remark that, if $a\in D_W$ we have $\partial_i(af_{S,W}/\prod_{(a,\chi)\in S}(1-a\chi) )=a\partial_i f_{S,W}/\prod_{(a,\chi)\in S}(1-a\chi) )$ if $i\leq k$ thus $a\partial_i\in D_W$. Instead, if $i>k$  we have  $\partial_i(a)=[\partial_i,a]\in D_W$ and  $\partial_i f_{S,W}/\prod_{(a,\chi)\in S}(1-a\chi) )=0$ and still the claim follows.

 \begin{lemma}

     For each component $W$, the complex $M(W)\otimes\Xi$ is isomorphic
      to the tensor product of two complexes $M_1$ and $M_2$.

       $M_1$ is the De Rham complex for the variety $W$ which is isomorphic as an algebraic variety, to a $r-k$ dimensional torus.

\smallskip

$M_2= \mathbb C[\partial_1,\ldots ,\partial_k] \otimes \bigwedge(\omega_1,\dots,\omega_k)\otimes  V_{W},$  with differential
$$\delta(p\otimes\omega\otimes v)=\sum_{i=1}^k\partial_ip\otimes\omega_i\wedge\omega\otimes v$$
for $p\in \mathbb C[\partial_1,\ldots ,\partial_k]$, $\omega\in \bigwedge(\omega_1,\dots,\omega_k)$, $v\in \tilde V_{W}.$
      $$H^{s}(M_2)=\begin{cases}0 {\rm \ if\ } s<k\\  V_{W}\otimes \omega_1\wedge\cdots\wedge \omega_k{\rm \ if\ } s=k\end {cases}$$

      \end{lemma}

    \proof

 Using the choices made above, we can identify the coordinate ring of $W$ with $\mathbb C[\xi_{k+1}^{\pm 1},\ldots ,\xi_r^{\pm 1}]$. Thus $M_1:= \mathbb C[\xi_{k+1}^{\pm
1},\ldots ,\xi_r^{\pm 1}]\otimes \bigwedge(\omega_{k+1},\dots,\omega_r)$, with the usual De Rham differential, is the De Rham complex of $W$.

On the other hand,  take  $f\otimes \phi_1\in M_1,\  p\otimes\phi_2 \otimes v\in M_2$ (with the obvious notations). It is immediate to verify that the product  $$f\otimes \phi_1\otimes p\otimes\phi_2 \otimes v\mapsto fp\otimes \phi_1\wedge\phi_2 \otimes v$$
 induces an isomorphism of complexes $M(W)\otimes\Xi\simeq
M_1\otimes M_2$.

The complex $M_2$ is the tensor product of the vector space $V_{W} $ (thought of as trivial complex in degree 0) and the complexes  concentrated in degrees 0,1, $\mathbb C[\partial_i]\to \mathbb C[\partial_i] \omega_i$ with $d(p)=\delta_ip\omega_i$,  whose cohomology is generated by the class of $\omega_i$ in degree 1. The statement about the cohomology of $M_2$ follows immediately from this.\qed

As an immediate consequence we obtain in a more explicit form the result of  Looijenga, \cite{lo}, \S 2.4.3.
\begin{theorem}\label{decomp}
   For each integer $i\geq 0$, we have a (non canonical) decomposition, as $W$ runs over the components of the arrangement:
    $$H^i(\mathcal A_{\Delta})=\oplus_{W} H^{i-{\rm codim }W}(W)
    \otimes  V_{W}.$$ \end{theorem}
 More precisely a basis for the cohomology  can be  given by the classes of the forms
\begin{equation}\label{decomp10}{f_{S,W}\over (1-a_1\chi_1) \cdots   (1-a_k\chi_k)}d\log\xi_1\wedge\dots\wedge d\log\xi_k\wedge d\log \xi_{j_1}\wedge\cdots\wedge
d\log \xi_{j_s}.\end{equation}
Here  $W$ runs over the set of components,  $S=\{(a_1,\chi_1)<\cdots <(a_k,\chi_k)\}$  runs over the non broken circuits in $\Delta_{W}$, and for given $(S,W)$,
$f_{S,W}\in \mathbb C[\Sigma ]$ is a chosen  representative of the  primitive idempotent $\phi_{S,W}$. Finally  $\xi_1,\dots,\xi_r$ is a basis of the character group such that the $\chi_i$ are in the subgroup spanned by the first $k$ elements.

\begin{remark}\label{base}   1) The previous forms do not depend  on the choice of a complement of the lattice spanned by $\xi_1,\dots,\xi_k$ in $\Lambda$ but only on the choice of a basis of $\Lambda/\langle\xi_1,\dots,\xi_k\rangle$.\smallskip

2)The filtration of $R$ by the subspaces $R_k$, induces for each $i$, a natural filtration  $H^i(\mathcal A_\Delta)$  on cohomology such the $k^{th}$ graded component is. 
 $$H^i(\mathcal A_{\Delta})_k/H^i(\mathcal A_{\Delta})_{k-1}=\oplus_{W |\, {\rm codim } \, W=k}\  H^{i-k}(W)
    \otimes  V_{W}.$$
    
3) Using Theorem \ref{decomp} one can write down an explicit formula for the Poincar\'e polynomial of $\mathcal A$ in terms of the poset of components. It is easily seen that this formula coincides with the formula given by Looijenga in \cite {lo}, 2.4.3.
\end{remark}
 With the previous notations we can also introduce a different basis which will be useful for the future computations:

\begin{proposition}\label{leforme}  The forms
$$f_{S,W}d\log(1-a_1\chi_1)\wedge\cdots\wedge d\log(1-a_k\chi_k)\wedge d\log \xi_{j_1}\wedge\cdots\wedge
d\log \xi_{j_s}$$  (where $k<j_1<\dots<j_s\leq r$), are closed and their classes give a basis of the cohomology of $\mathcal A_{\Delta}$.
\end{proposition}\proof The fact that these forms are closed follows as in the previous arguments.  These forms are obtained by an invertible triangular change of basis with respect to the forms (\ref{decomp10}).  In fact 
$$f_{S,W}d\log(1-a_1\chi_1)\wedge\cdots\wedge d\log(1-a_k\chi_k)= f_{S,W}\prod_i{-a_i \chi_i\over (1-a_i\chi_i)  }d\log \chi_i.$$
Substitute for $-a_i \chi_i$ the sum $(1-a_i \chi_i)-1$ and develop all the terms getting the required triangular expression.\qed

For our work the top dimension $i=r$ is particularly important. In this case the forms 
\begin{equation}\label{tfor}f_{S,W}d\log(1-a_1\chi_1)\wedge\cdots\wedge d\log(1-a_k\chi_k)\wedge d\log \xi_{j_1}\wedge\cdots\wedge
d\log \xi_{j_{r-k}}\end{equation} are independent of the way in which we complete the basis  $\xi_1,\dots,\xi_k$ except for the orientation which we fix in an arbitrary way. Their classes will be denoted by $\omega_{S,W}.$ Such a class is said to be {\it associated} to $W$.

Given a top cohomology class $\psi$ when we expand $\psi$ in the basis $\omega_{S,W} $ we call the coefficients the {\it local residues} of $\psi$ at $S,W$ and write:
$$\psi=\sum_{S,W}res_{S,W}(\psi )\omega_{S,W}. $$

In particular, the form \begin{equation}\label{can}\omega_{\emptyset,T}=d\log \xi_1\wedge \dots \wedge d\log \xi_r\end{equation} depends only on the orientation of the basis $\xi_1,\dots,\xi_r$ of $\Lambda$ and will be denoted by $\omega_T$ (or just $\omega$).

 \subsection{The local residue}

 Let $P$ be one of the points of the arrangement,  thus we have  as usual a list $\Delta_P$ of the divisors $1-a_i\psi_i$ passing through $P$.
 
 Consider  $\mathfrak p$ the tangent space of $T$ in 1 and let $e_P:t\mapsto Pexp(t)$  which gives local coordinates around $P$.

 The preimage of one of the divisors $1-a_i\psi_i=0$ is a hyperplane whose equation $h_i$ is defined by  $\psi_i(e^t)=e^{h_i(t)}$.  We can restrict to a small enough ball of center 0, so that the preimage of $\mathcal A_\Delta$ is exactly the complement of the intersection of these hyperplanes with this ball. 
 
 The total ordering on $\Delta_P$ transfers to a total ordering of the vectors $h_i$ which preserves the notion of no broken circuits. From a topological point of view this is diffeomorphic to the complement of this hyperplane arrangement. 
Let us now compute  now the pull back under $e_P$ of the top forms
and top cohomology classes $\omega_{S,W}$ defined in  (\ref{tfor}). 
 \begin{theorem}\label{primc}
1)  Under $e_P^*$ all cohomology classes which are not associated to $P$ go to 0.

2)    The cohomology class of the form  
 $$e_P^*(f_{S,P}d\log (1-a_1\psi_1)\wedge\ldots\wedge(1-a_r\psi_r))$$ equals that of the form 
 $d\log h_1\wedge\ldots\wedge d\log h_r$.
 \end{theorem} 
 \proof Consider one of the forms $d\log(1-a\psi)$  if $(1- a\psi )(P)\neq 0$ we can choose (locally) a determination of the logarithm so if the ball is small enough the pull back  of this form is exact.  The same applies to the forms $d\log\chi$.   If instead $(1-a\psi )(P)= 0$, with $\psi=e^h$  we have that:
 $$e_P^*d\log (1-a\psi ) =d\log h+\gamma, $$
 with $\gamma$ holomorphic on the ball.
 
 One concludes with two remarks.  
 
 A product  of $<r$ factors of type $d\log h_i$ times a holomorphic  $r-k$ form is exact. 
 
 The form   $f(t)d\log h_1\wedge\ldots\wedge d\log h_r$ is cohomologous to $f(0)d\log h_1\wedge\ldots\wedge d\log h_r$.  
  \qed
 
It is well known \cite{te}, that the forms $w_S:= d\log  h_1\wedge\ldots\wedge d\log h_r$ as $S:=\{h_1,\ldots,h_r\}$ runs over the no broken circuits of the hyperplane arrangement, form a basis of the cohomology of the complement of the hyperplanes. We get

   \begin{theorem}\label{loc}  Given a point $P$ of the arrangement and a differential form of top degree $\Psi$ we have  
   
\begin{equation}\label{clocal}e_P^*([\Psi])=\sum_{S\in\Delta_P}res_{S,P}([\Psi] )w_{S}. \end{equation}
where $[\Psi]$ denotes the cohomology class of $\Psi$.

 \end{theorem}

 In  \cite{dp1}  we have constructed   geometrically a dual homology basis, and shown how one can explicitly compute the integral of the form on each element of the basis as a local residue around a new point at {\it infinity}  in a  suitable blow up  which transform the divisor of hyperplanes into one with normal crossing.
  This then becomes a method to compute the local residues $res_{S,P}([\Psi] )$.

  In fact we will need to analyze, for suitable $\Psi$,  the function $$c(\beta):=res_{S,P}([\beta^{-1}\Psi ])$$ as $\beta\in\Lambda$. Let us briefly discuss an important qualitative feature of this number as a function of $\beta$.
  
  \begin{definition}
 A function $g $,   on $\Lambda$   will be called  a {\bf periodic polynomial} if it is  a polynomial   on each coset of some sublattice   $\Lambda'$    of  $ \Lambda$.
\end{definition}

Remark that as sublattice we can take   $n\Lambda$ for some positive integer $n$.
If $f$ is a periodic polynomial on  $ \Lambda$, it is also a periodic polynomial on  $ 1/n\Lambda$ for every positive integer $n$.
 \smallskip

\begin{proposition}\label{periodic} 
$c(\beta)=res_{S,P}([\beta^{-1}\Psi ])$ is a periodic polynomial on $ \Lambda$.
\end{proposition}
\proof  In the coordinates $y_1,\ldots ,y_r$ of the blow up (see \cite{dp1} for the explicit formulas), the form 
$\Psi$ develops as: 
$$\Psi={1\over y_1^{h_1}\cdots y_r^{h_r}}G(y_1,\ldots ,y_r)dy_1\wedge\cdots\wedge dy_r$$
where $G$ is a holomorphic function near zero. As for $\beta^{-1}$ it can be written as $\beta^{-1}(p)$ times a convergent power series in the $y$'s with coefficient which depend polynomially on $\beta$. From this the claim follows since $\beta^{-1}(p)$ is constant on the sublattice $\Lambda'$ spanned by $\psi_1,\ldots ,\psi_r$.\qed

  In the way we have organized our discussion,  the computation of the local residues is performed by viewing the toric arrangement  locally, around each of the points $P$  of the arrangement. We have used the fact that locally the arrangement  is isomorphic to a hyperplane arrangement, to which we then apply our geometric theory of blow ups. In fact, although we do not use this explicitly,  the theory is even global. It is clear that the stratification of a torus induced by the toric arrangement is {\it conical} as defined in \cite{mp}. The theory of irreducible strata and nested sets can be computed locally and thus reduced to the hyperplane case where it is well understood. Then, by applying  the theory developed in \cite{mp},
we can construct a smooth model  $\tilde T$ with a proper map $\pi:\tilde T\to T$. $\pi$  is isomorphic on the complement $\mathcal A_{\Delta}$ of the arrangement. The divisor in  $\tilde T$,   complement of $\mathcal A_{\Delta}$,  has normal crossings and it has globally the combinatorics described locally, for all the points $P$ of the arrangement, by the corresponding hyperplane arrangements.

\subsection{Coverings} Let us now consider a finite $n-$sheeted covering $\pi: U\to T$ of tori. It corresponds to an embedding $\Lambda\subset M$  of character groups, with $\Lambda$ of index $n$ in $M$.  It is a Galois covering with Galois group the dual $(M/\Lambda)^*$ of $M/\Lambda$.

The set $\Delta$ clearly defines also an arrangement in $U$,  to every  component $W$ of the arrangement in $T$ correspond the connected components of $\pi^{-1}(W)$. Let us denote by $\mathcal A_\Delta^T, \mathcal A_\Delta^U$ the complement for the two arrangements, we have thus a map $\pi:\mathcal A_\Delta^U\to
 \mathcal A_\Delta^T$ and in cohomology $\pi^*:H^*(\mathcal A_\Delta^T)\to H^*( \mathcal A_\Delta^U).$\smallskip
 
 The vector spaces of invariant derivations for $T$ and $U$ coincide. Thus we can take a common basis $\partial_1,\dots,\partial_r$ and we have \begin{equation}D_U=\mathbb C[U][\partial_1,\dots,\partial_r]=\mathbb C[U]\otimes_{\mathbb C[T]}D_T.\end{equation}
For a $D_T$ module $M$ we can construct the {\it induced module} on $U$  as:
$$D_U\otimes_{D_T}M= \mathbb C[U]\otimes_{\mathbb C[T]}M. $$
The action of a derivation $\partial_i$ on   $a\otimes m, \ a\in \mathbb C[U]$ is clearly:
$$\partial_i(a\otimes m)=\partial_i(a)\otimes m+a\otimes \partial_i(m). $$
The following Proposition, describing the induction of the modules  $N(W)$ of \S 2.2,  is only a reformulation of  Lemma \ref{induc}:

\begin{proposition}  Given  a covering $\pi:U\to T$ and   an irreducible component $W\subset T$ of an arrangement, there is a canonical isomorphism of $D_U$ modules
$$ \mathbb C[U]\otimes_{\mathbb C[T]}N(W)\simeq \oplus_{i=1}^hN(W_i).$$

Where the $W_i$ are the irreducible components of $\pi^{-1}(W)$.
\end{proposition}

This has an important consequence for the canonical filtration of $  \mathbb C[U]_\Delta:=  \mathbb C[U][d^{-1}]$.
Start from the fact that $  \mathbb C[U]_\Delta=\mathbb C[U]\otimes_{\mathbb C[T]} \mathbb C[T]_\Delta$. 
\begin{theorem}   For each $k=0,\ldots r$:
$$  (\mathbb C[U]_\Delta)_k= \mathbb C[U]\otimes_{\mathbb C[T]} (\mathbb C[T]_\Delta)_k.$$

$$  (\mathbb C[U]_\Delta)_k/(\mathbb C[U]_\Delta)_{k-1}= \mathbb C[U]\otimes_{\mathbb C[T]} [(\mathbb C[T]_\Delta)_k/(\mathbb C[T]_\Delta)_{k-1}] $$
  hence we have   canonical isomorphisms:  
$$ (\mathbb C[U]_\Delta)_k/(\mathbb C[U]_\Delta)_{k-1}\simeq \oplus_{W} \mathbb C[U]\otimes_{\mathbb C[T]} N_T(W)\otimes V_{W}\simeq \oplus_{W}\oplus_Z N_U(Z)\otimes V_{W}$$
$W$ runs over the components of the arrangement in $T$ of codimension $k$ and, for given $W$, $Z$ runs  over the components of the arrangement in $U$ which lie over $W$.\end{theorem}
\proof
Since $\mathbb C[U]$ is a free module over $\mathbb C[T] $ the functor $\mathbb C[U]\otimes_{\mathbb C[T]}-$ is exact hence the claims follow since the $k-$th level of the filtration of $\mathbb C[U]_\Delta$ can be characterized as the maximal $D_U$ submodule whose characteristic variety is the union of conormal bundles of subvarieties of codimension $\leq k$.
\qed

      Since $\Lambda_{\mathbb R}=M_{\mathbb R}$ we fix once and   for all an  orientation of this vector space which induces compatible orientations for all coverings.   The we have the canonical classes $\omega_T,\omega_U$ as in formula  \ref{can} and:
   \begin{lemma}\label{multipl}  If $\pi$ is  of degree $n$: $$\pi^*(\omega_T)= n\omega_U$$
   \end{lemma}
   \proof
   By the theory of elementary divisors, there is an   oriented basis  $ \mu_1,\dots,\mu_r$ of $M$ and positive integers $n_1,\dots,n_r$ such that  $n=\prod_in_i$ and   $\mu_1^{n_1},\dots, \mu_r^{n_r}$ is a basis of $\Lambda$. Using these two bases  the claim is clear. 
   \qed  
   
 Let us compute in general $\pi^* (\omega_{S,W})$ where $\omega_{S,W}$ is defined in formula (\ref{tfor}). 
 
 Given the lattice $\langle S\rangle\subset\Lambda$ its associated $\overline{\langle S\rangle}$ of elements of $\Lambda$ which are torsion modulo $\langle S\rangle$ and finally $ \tilde{\langle S\rangle}$ of elements of $M$ which are torsion modulo $\langle S\rangle$ we can construct summands:
 $$\Lambda= \overline{\langle S\rangle}\oplus\Pi,\quad M= \tilde{\langle S\rangle}\oplus \tilde\Pi,\quad \Pi\subset  \tilde\Pi$$
We have that $\pi^{-1}(W)$ has $[ \tilde{\langle S\rangle}:\overline{\langle S\rangle}]$  components while, if  $\xi_1,\dots,\xi_{r-k}$ resp. $\eta_1,\dots,\eta_{r-k}$ are bases of $ \Pi$ resp. of $  \tilde\Pi$ with the same orientation, we have 
\begin{equation}\label{ciccio}d\log\eta_1\wedge \dots \wedge  d\log \eta_{r-k}=p^{-1}d\log\xi_1\wedge \dots \wedge d\log \xi_{r-k},\quad p=[ \tilde\Pi:\Pi].\end{equation}
For every component $W'$ of $ \pi^{-1}(W)$ we can clearly choose a representative, $f_{S,W'}$ of the corresponding primitive idempotent  in such a way that
$$f_{S,W}=\sum_{W'\text{\ component of\ \ } \pi^{-1}(W)}f_{S,W'}.$$
We then set 
$$\omega_{S,W'}:=f_{S,W'}d\log(1-a_1\chi_1)\wedge\cdots\wedge d\log(1-a_k\chi_k)\wedge d\log \eta_{1}\wedge\cdots\wedge d\log \eta_{r-k}$$

\begin{proposition}
We have that:
$$ \pi^* (\omega_{S,W})=p^{-1}\sum_{W'}\omega_{S,W'}$$
as $W'$ runs over the  components of $ \pi^{-1}(W)$.
\end{proposition}
\proof
Everything follows from  identity (\ref{ciccio}).
\qed

Notice that in the case in which $W$ is a point we have that $p=1$.

  Let us specialize to one important case. Let $U$ be the torus with character group  $M={1\over m}\Lambda$ for some integer $m$  so that  the degree of the map $\pi:U\to T$ is $m^r$.

   Inverting  the elements of $\Delta$ is equivalent to inverting the elements of:
   $$\Delta^{1/m}:=\{(b,\sigma ),\ b\in\mathbb C^*,\  \sigma \in  {1\over m}\Lambda ,\ |\, (b^m,\sigma^m)\in \Delta\}.$$
    In fact let $(a,\chi)\in \Delta$, choose $b$ with $b^m=a$. We have $m$ pairs $(\zeta^ib,\sigma),\ i=0,\dots,m-1,\ \zeta=e^{2\pi i/m} $ with $((\zeta^ib)^m,\sigma^m)=(a,\chi)$   where $\sigma$ is uniquely determined.   Moreover
    $$1- a\chi=\prod_{ i=0,\dots,m-1}(1-b\zeta^i\sigma).$$
Take  a component $W\subset T$ of the arrangement $\mathcal R_\Delta$.   This is  a  connected component of the variety  of   equations   $a_i\chi_i=1\ i=1,\dots,k$, where $(a_i,\chi_i)\in\Delta$ and the $\chi_i$ are  linearly independent.

 $\pi^{-1}(W)$  consists  of $m^k$ connected components, one component for each of the varieties $\zeta^jb_i\sigma_i=1,\ i=1,\dots,k$ as $b_i^m=a_i,\sigma^m=\chi_i$ and the exponents $j$ vary in the $m^k$  possible different ways.

Recall that $\overline\Delta$ has been defined as the set of characters $\chi$ such that there is $a\in\mathbb C^*$ with $(a,\chi)\in\Delta$ and define $\overline\Delta^{1/m}$ in the same way. Notice that $\overline\Delta^{1/m}$ is in canonical bijection with $\overline\Delta.$  Thus the chosen ordering on  $\overline\Delta $ gives an ordering on $\overline\Delta^{1/m}$.

 For the given component $W$,     select  a component $W'$  of $\pi^{-1}(W)$. 

  It is then immediate to see  that the map $(b,\sigma)\mapsto (b^m,\sigma^m)$ defines a bijection between $\Delta^{1/m}_{W'}$ and $\Delta_W$.     This establishes also a bijection between the  no broken circuits on $W'$  and the ones on  $W$. 
  
   As for cohomology, take a basic form
   $$f_{S,W}d\log(1-a_1\chi_1)\wedge\cdots\wedge d\log(1-a_k\chi_k)\wedge d\log \xi_{j_1}\wedge\cdots\wedge
d\log \xi_{j_s}.$$  Choose as   complement of the lattice associated   the elements $\chi_i$ or their $m-$th roots,    in the character lattice of  the covering, the one  generated by the $m-$ roots $\eta_i$ of the $\xi_i$.

Since $$d\log(1-a_i\chi_i)=d\log(\prod_{j=0}^{m-1}(1-\zeta^jb_i\sigma_i)=\sum_{j=0}^{m-1}d\log(\ 1-\zeta^jb_i\sigma_i),$$ $$ d\log\xi_i=d\log(\eta_i^m)=md\log\eta_i,$$ we have that the previous form is a sum of terms
\begin{equation}\label{ess}m^{s}f_{S,W}d\log(1-\zeta^{j_1}b_1\sigma_1)\wedge\cdots\wedge d\log(1-\zeta^{j_k}b_k\sigma_k)\wedge d\log \eta_{j_1}\wedge\cdots\wedge
d\log \eta_{j_s}\end{equation}
The function $f_{S,W}$ has value 1 on the preimage in $U$ of $W$ and value 0  on the preimage of all the other components of the variety of equations  $1-a_i\chi_i=0$.

Thus, with the previous notations, for $S_{\underline j}:= \{ 1-\zeta^{j_1}b_1\sigma_1 ,\cdots, 1-\zeta^{j_k}b_k\sigma_k \}$ we have that:
\begin{lemma}  We can choose representatives $f_{ S_{\underline j}}$
 as $W_t$ runs over all the components of the variety of equations $1-\zeta^{j_1}b_1\sigma_1 =\cdots= 1-\zeta^{j_k}b_k\sigma_k=0$ in $U$, which lie over $W$ so that:
 $$f_{S,W}=\sum_t f_{ S_{\underline j} ,  W_t} $$
 \end{lemma}
 \proof
 We choose all except one which we define from the previous formula and verify that it is indeed a required representative.
 \qed
 
Remark we should compare the two expressions associated to a component in $U$ using either $\Delta$ or $\Delta^{1/m}$.

We have 
$$ d\log(1-\zeta^{j_1}b_1\sigma_1)\wedge\cdots\wedge d\log(1-\zeta^{j_k}b_k\sigma_k)=$$$$\prod_{i=1}^k\prod_{j=0, j\neq j_i}^{m-1}(1-\zeta^{j }b_k\sigma_k)d\log(1-a_1\psi_1)\wedge\cdots\wedge d\log(1-a_k\psi_k)$$

A special case will be needed, when $W=P $ reduces to a point. The 0-dimensional components of the arrangement play a special role in the final formulas. For such $P $  in formula (\ref{ess}) we have $s=0$   and thus:
  \begin{equation}\label{esss} f_{S,P }d\log(1-a_1\chi_1)\wedge\cdots\wedge d\log(1-a_r\chi_r)=\end{equation}
  $$\sum_{Q \in\pi^{-1}(P )} f_{ S_{\underline j} ,  Q }d\log(1-\zeta^{j_1}b_1\sigma_1)\wedge\cdots\wedge d\log(1-\zeta^{j_r}b_r\sigma_r).$$
In the right hand side of this formula  the   components lying over $W=P $ are just the $m^r$ points in $Q \in\pi^{-1}(P )$  the fiber.  Each such point $Q$ determines a unique    no broken circuit  $S_{\underline j}$ obtained as  $m-$th root  of $S$.  The formula induces a corresponding formula in cohomology.
\section{The cohomology algebra}
\subsection{A basic identity}
We start from a basic formal identity:
\begin{equation}\label{iden} 1-\prod_{i=1}^nx_i=\sum_{I \subsetneq \{1,2,\dots,n\}}\prod_{i\in I}x_i\prod_{j\notin I}(1-x_j).\end{equation}
The proof  is by induction on $n$.

We split the sum in 3 terms: $I=\{1,\dots,n-1\},\ I\subsetneq \{1,\dots,n-1\}$ and finally $n\in I$. We get
$$ \prod_{i=1}^{n-1}x_i(1-x_n)+(1- \prod_{i=1}^{n-1}x_i)(1-x_n)+x_n( 1-\prod_{i=1}^{n-1}x_i)=1-\prod_{i=1}^nx_i.$$
Formula (\ref {iden})  implies

 $$   \sum_{I \subsetneq \{1,2,\dots,n\}} {1\over (1-\prod_{i=1}^nx_i)}\prod_{i\in I} {x_i\over  (1-x_i)}={1\over \prod_{i=1}^n(1-x_i)}.$$

 We want to interpret this formula as an identity between certain differential forms.

 Set, for $i=1,\ldots ,n$,  $\omega_i:=d\log(1-x_i),\ \psi_i:=d\log x_i$,. Also for $s=0,\ldots n$, set  $\theta^{(s)}=d\log (1- \prod_{i=1}^sx_i^{-1}\prod_{j=s+1}^nx_j)$.  If we take a  proper subset $I=\{i_1<\cdots <i_t\}$ in $\{1,\ldots ,n\}$ and let $J=\{j_1<\cdots <i_{n-t}\}$ be its complement, we can then define the $n$-forms:
$$\Phi_I^{(s)}=(-1)^{s_I}\omega_{i_1}\wedge\cdots\wedge \omega_{i_t}\wedge \psi_{j_1}\wedge\cdots\wedge \psi_{j_{n-t-1}}\wedge \theta^{(s)}$$
and
$$\Psi_I=(-1)^{s_I}\omega_{i_1}\wedge\cdots\wedge \omega_{i_t}\wedge \psi_{j_1}\wedge\cdots\wedge \psi_{j_{n-t}}$$
with $s_I$ equal to the parity of the permutation $(i_1,\ldots ,i_t, j_1,\ldots ,i_{n-t})$.

We have
\begin{proposition}\label{dritto}
For each $0\leq s\leq n$,  the $n$-form
$\omega_1\wedge\cdots\wedge\omega_n$ can be written as a linear combination with integer coefficients of the $n$-forms $\Phi_I^{(s)}$ and $\Psi_I$.
\end{proposition}
\proof We first deal with the case $s=0$. In this case, computing we get
 \begin{equation}\Phi_I^{(0)}=
(-1)^{|I|} \prod_{i\in I} {x_i\over  (1-x_i)}{1\over (1-\prod_{i=1}^nx_i)}dx_1\wedge
 dx_2\wedge\dots dx_n.\end{equation}
 Thus our identity (\ref {iden}) can be translated into: \begin{equation}\label{casozero} \sum_{I \subsetneq
 \{1,2,\dots,n\}}(-1)^{|I|+n}\Phi_I^{(0)}=\omega_1\wedge\cdots\wedge\omega_n.\end{equation}
 proving our claim.

 In the general case let us observe  that $$d\log(1-x^{-1})=d\log{x-1\over x}=-d\log (x)-d\log (1-x).$$
   Therefore, the substitution of  i$x_i$ with $x_i^{-1}$ for $i=1,\ldots ,s$, corresponds in   the formula (\ref{casozero}) to substitute $\omega_i$ with $-\omega_i-\psi_i$ getting a new formula which we can denote as (\ref{casozero}'). \qed

\subsection{Formality}\label{unim}   In this section we shall assume that our set  $\Delta$ is {\it unimodular}, with the terminology and notations of \S  3.1, this means that   given a  subset  $\Psi\subset \Delta$ which is linearly independent, $\pi(\Psi)$ generates a lattice which is a
direct factor in $\Lambda$. In particular, if $\xi_0$ is dependent on $\xi_1,\dots ,\xi_k$  with all these elements in $\Delta$ we must have that $\xi_0$ is a product of the $\xi_i$ with exponents $\pm 1$.   We call this a {\it simple dependency relation}.   More generally we say that the element
 $(a_0,\xi_0)$ is dependent on $(a_1,\xi_1),\dots, (a_k, \xi_k)$ if $\xi_0$ is a monomial in the $\xi_i,\ i=1,\dots,k$ and $a_0$ is the value of the same monomial in the $a_i$.
\smallskip

Let us start with a formal construction. Define an exterior algebra in the following generators:

 A generator $\lambda_{(a,\chi)}$ for every $(a,\chi)\in\Delta$ and a generator $\mu_\chi$ for every character $\chi\in\Lambda$.  
 \smallskip
 
 We then impose a set of relations.  Of this first we have the essentially trivial relations:
 \begin{enumerate}
 \item 
$\mu_{\chi_1\chi_2}=\mu_{\chi_1}+\mu_{\chi_2}$, for all $\chi_1,\chi_2\in\Lambda$.
\item  If $\chi_1,\dots,\chi_s$ are dependent the product of elements $\sigma_i$ where $\sigma_i=\mu_{\chi_i}$ or $\sigma_i=\lambda_{(a_i,\chi_i)}$ is 0.  
\end{enumerate}   We come now to the main relations, which generalize the ones for hyperplanes.

Given $(a_0,\chi_0),(a_1,\chi_1),\dots,(a_n,\chi_n)\in\Delta$  with$$ a_0\chi_0=\prod_{i=1}^s(a_i\chi_i)^{-1}\prod_{j=s+1}^na_j\chi_j.$$ 
We take the formula (\ref{casozero}')  and substitute $\psi_i$ with $\mu_{\chi_i}$ and $\omega_i$ with $\lambda_{a_i,\chi_i}$. We get then a formal expression which we impose as a new relation. 

We do this for all the simple dependency relation and we call $\mathcal H$ the resulting quotient algebra.
\smallskip

Consider the subalgebra $H$ in $\Omega=R\otimes\Xi$ generated by the 1-forms $d\log(1-\nobreak  a\chi), d\log \chi$. $H$ clearly consists of closed forms so that we obtain a algebra homomorphism
$$f:H\to H^*(\mathcal A_{\Delta}).$$
 It is also easy to verify, using Proposition \ref{dritto}, that we have  a homomorphism  $g:\mathcal H\to H$ given by:
 $$ g(\lambda_{(a,\chi)}):=d\log(1- a \chi ) ,\ g(\mu_\chi):=d\log\chi.$$

 \begin{theorem}\label{formality}   The homomorphisms $g,f$ are isomorphisms.
 \end{theorem}
  \proof The assumption that $\Delta$ is unimodular clearly implies that for each component $W$ and for any $S$ associated to $W$, the defining ideal of $W$ is generated by the elements $1-a\chi$ with $(a,\chi)\in S$. In particular  $S$ is associated only to one component and thus the idempotent $\phi_{S,W}=1$. From this and Remark \ref{base} we deduce that $f$ is surjective.
  
  The fact that $g$ is surjective is clear from its definition. So, in order to prove our claim, it suffices to see that $fg$ is injective.

 We fix a total ordering on $\Delta$ in order to apply the theory of no broken circuits.
The first set of  relations implies that the subalgebra generated by the elements $\mu_{\chi}$ is a homomorphic image of the exterior algebra $\Xi$. So we can consider $\mathcal H$ as a $\Xi$-module.

For a fixed component $W$ we define the subspace $\mathcal V_{W}\subset \mathcal H$ spanned by the monomials $\lambda_S=\prod_{(a,\chi)\in S}\lambda_{(a,\chi)}$ (here the product is taken according to the fixed ordering of $\Delta$)  as $S$ runs among the subsets of  $\Delta$
 associated to
$W$. 
The first set of  relations also implies that $\mathcal V_{W}$ is annihilated by the elements $\mu_{\chi}$ if the character $\chi$ is constant on $W$.

Given a monomial $m$ in the generators $\lambda_{(a,\chi)}$ and $\mu_{xi}$, we  define its weight as the number of factors of $m$ of the first type    $\lambda_{(a,\chi)}$.  Notice that the weight of $\lambda_S\in \mathcal V_{W}$ equals the codimension of $W$. Then, whenever $S$ is  a broken circuit, the relations of the second type   allow to replace it by a product of elements with lower weight or lower in the lexicographical order. 
This implies that modulo elements of smaller weight any element in $\mathcal V_W$ can be written as a linear combination of elements $\lambda_S$ with $S$  a non broken circuit associated to  $W$. 
  
 From this our claim follows immediately from Proposition \ref{leforme}.
 \qed

  \begin{remark} We wish to point out that our result shows  in particular that in $H^*(\mathcal A_{\Delta})$  the algebraic relations between the generating forms $d\log (1-a\chi)$ and $d\log\xi$ resemble, but are more complicated, than the relations of Orlik-Solomon in the case of hyperplane arrangements.
\end{remark}
  \section{Residues}

The rest of this paper is devoted  to show the relationship between our work and the theory of counting integer points in polytopes developed by Brion, Szenes, Vergne (\cite{BV1},\cite{SV1}).  In a way, it is mostly a  reformulation in our language of their results and does not claim to be particularly original.\smallskip

Our main point is that we want to give a geometric interpretation  of the sum over certain roots of 1 which appears in their formulas. In our approach this has a clear geometric explanation, first as sum over the zero dimensional components of the toric arrangement and finally as a further sum over suitable {\it points at infinity} in the blown up model where the arrangement is made of divisors with normal crossing. In this sense we recover a {\it residue formula} as in the 1-dimensional case.  The points at infinity which contribute with a non zero residue depend on the character we choose (in the partition function) in a way which is encoded in the notion of Jeffrey--Kirwan residue.\smallskip

    From now on we shall assume that there is  a one parameter group $j:G_m\to T$ such that for each $(a,\gamma)\in\Delta$, $\langle x_i\,|\, j(t)\rangle=t^{m_i}$ with $m_i>0$.
    
    This assumption implies   the following geometric considerations. 
    
    Let $\Lambda_{\mathbb  R}:=\Lambda\otimes\mathbb  R$, the real  the vector spanned by the characters. In  $\Lambda_{\mathbb  R}$  denote by  $C(\Delta)$, the convex cone of non negative linear combinations of the elements $\gamma$ for all  $(a,\gamma)\in\Delta$.  The previous hypothesis implies that $C(\Delta)$ is a {\it pointed cone}, that is it lies entirely one one side of some  hyperplane.

    We can thus   consider the ring $P$ of formal power series     in the characters $\chi\in C(\Delta)$  and  its localization $Q$  obtained inverting all characters.

    Let us denote the De Rham complex for $Q$  by $\Omega^*(Q)$.
    Using the form $\omega:=\omega_T=d\log\xi_1\wedge\dots \wedge d\log \xi_r$ where $\xi_i$ is an oriented basis of $\Lambda$, we have that the top forms  can be  identified to $Q\omega$.
   \begin{lemma} $H^r(\Omega^*(Q))=\mathbb C$ and it is generated by the class of $\omega$. \end{lemma}
    \proof  Notice that,  setting as before  $\partial_i={  \xi_i }{\partial\over \partial_{\xi_i}}   $, we have that, if $\chi=\prod_i\xi_i^{n_i}$, $\partial_i\chi=n_i\xi$ so that the characters are eigenvectors for the $\partial_i$. The differential has the  form $\psi\mapsto \sum_i\partial_i\psi\wedge  d\log\xi_i$  we claim that the exact forms  of degree $r$ are exactly the ones which can be written as $f\omega$ where $f$ is a series in the characters without constant term. In one way it is clear, a form $\sum_i\partial_i\psi\wedge  d\log\xi_i$ has no constant term, conversely let $f\in Q$ write $f=f_1+ g$ where $f_1=\sum_{h=-N,\ h\neq 0}^\infty  \xi_1^{h }a_h $ is the sum of all terms containing $\xi_1$ with a non 0 exponent. Then $$f_1\omega=d\sum _{h=-N,\ h\neq 0}^\infty  h^{-1}\xi_1^{h }a_h    d\log\xi_2\wedge  d\log\xi_3\ldots \wedge  d\log\xi_r $$
   We then continue in the same way by induction with $g$. From this everything follows.
   \qed

  \begin{definition} The canonical map  $\Omega^r(Q)\mapsto H^r(\Omega^*(Q))=\mathbb C$
   will be called the {\bf residue} and  denoted by $res$. 
   \end{definition}
   Notice that, given a function $f\in Q$  we have that $res(f):=res(f\omega)$ is the constant term of the Laurent series   $f$. We deduce that, given $f=\sum_{\chi\in\Lambda}a_\chi \chi$ in $Q$: \begin{equation}\label{residuo}a_\chi=res(\chi^{-1}f)=res(\chi^{-1}f\omega).\end{equation}
    We have also a canonical map  $\Omega^r(R)\mapsto H^r(\Omega^*(R)) =H^r(\mathcal A_\Delta)$.
This will be called the {\it total residue} and denoted    by $Tres$. Given a form $f\omega\in \Omega^r(R)$ sometimes we shall write just $Tres(f):=Tres(f\omega)$.

 Setting
 $$i({1\over 1-a\chi})=\sum_{n\geq 0}(a\chi)^n$$
 for each $(a,\chi)\in \Delta$, we get an inclusion morphism $i:R\mapsto Q$.  
 
In the next paragraph it will convenient to write the lattice $\Lambda$ additively and write the characters as  $e^\alpha$. Then we can interpret the map  $i:R\mapsto Q$  as a map from $R$ to a space of functions on     $\Lambda$ and call it the {\it  Discrete Inverse Laplace Transform}, or DILT   for short.
\smallskip

The map $i$  extends
 to  a morphism of De Rham complexes   $\Omega^*(R)\mapsto\Omega^*(Q)$ and thus of cohomology groups, which by abuse of notation, we shall still denote by $i$. In top cohomology this gives the basic formula $$res(i(f))=i(Tres(f)).$$
 
\subsection{Coverings} 
   We pass now to analyze this picture for a given covering $\pi:S\to T$, of degree $n$,  with Galois group  $G=ker(\pi)$. We then define the algebras, $R,P,Q$ for the corresponding tori and denote them by $R_S,R_T,P_S,P_T,Q_S,Q_T$.   The covering corresponds to an inclusion $\Lambda\subset M$ of character groups,. Since $\Lambda_{\mathbb R}=M_{\mathbb R}$ we fix once and for all an  orientation of this vector space which induces compatible orientations for all coverings.   Recall the basic formula $\pi^*(\omega_T)= n\omega_S$ for   the canonical classes $\omega_T,\omega_S$.   
   
   We then have   commutative diagrams:

   \[\begin{CD}   R_T @>i>>Q_T \qquad@.  @.    \Omega^r(R_T)@>Tres >> H^r(\Omega^*(R_T))@>i>> H^r(\Omega^*(Q_T) \\@V\pi^*VV\! \! \! \! \! \! \! \! \! \! \! \! \! \! @V\pi^*VV  ,\qquad @.   @V\pi^*VV @V\pi^*VV @V\pi^*VV \\
     R_S @>i>>Q_S \qquad@.@.  \Omega^r(R_S)@>Tres >> H^r(\Omega^*(R_S)) @>i>>H^r(\Omega^*(Q_S)
   \end{CD}\]

   Under identification   of $H^r(\Omega^*(Q_T)),H^r(\Omega^*(Q_S))$ with $\mathbb C$ we have the commutative diagram:

     \begin{equation}\label{diag}\begin{CD}     H^r(\Omega^*(Q_T)@=\mathbb C \\@V\pi^*VV  @VnVV \\ H^r(\Omega^*(Q_S)@=\mathbb C
   \end{CD}\end{equation}

  Let us go back to the special  important case of $S$   the torus with character group  $M={1\over m}\Lambda$ for some integer $m$  so that  the degree of the map $\pi:S\to T$ is $m^r$.
 
 Recall formula (\ref{esss}) defining top forms 
$$ f_{S,P}d\log(1-a_1\chi_1)\wedge\cdots\wedge d\log(1-a_r\chi_r)=$$
  $$\sum_{Q \in \pi^{-1}(P)} f_{ S_{\underline j} ,  Q }d\log(1-\zeta^{j_1}b_1\sigma_1)\wedge\cdots\wedge d\log(1-\zeta^{j_r}b_r\sigma_r).$$
 If  $\pi^*(S)$   denotes the set of no broken circuits   $S_{\underline j}=(\zeta^{h_1}b_1,\sigma_1),\dots , (\zeta^{h_r}b_r,\sigma_r),$ obtained as $m^{th}$ roots of $S$, we have seen that they are in 1-1 correspondence with the points of the fiber $\pi^{-1}(P)$.
We shall denote by  $\omega_{S,P}, \omega_{S_{\underline j},Q}$ the cohomology classes of these forms.  Thus we have that:
  \begin{equation}\label{pullback} \pi^*( \omega_{S,P})=\sum_{Q\in\pi^{-1}(P)}\omega_{S_{\underline j},Q}=\sum_{S_{\underline j}\in\pi^*(S)}\omega_{S_{\underline j},Q}.\end{equation}

\section{The counting formulas}

As we discussed in the introduction, one of the aims of this paper is to understand in a more algebraic geometric language, the counting formulas of \cite{BV}, \cite{BV1},  \cite{SV1} for integral points in polytopes.

In this setting we are given a family  $B$ of vectors in a lattice $\Lambda$  (which we treat additively) and, for $\beta\in\Lambda$ we want to count the number  $c_\beta$  of solutions of the equation:
$$\sum_{\alpha\in B}n_\alpha\alpha=\beta,\ \alpha\in\mathbb N.$$
This is clearly the value of a {\it partition function}, or in other words, the number of integral points of the polytope
$$\Pi_\beta:=\{(x_\alpha)\,| \sum_{\alpha\in B}x_\alpha\alpha=\beta,\ 0\leq x_\alpha\in\mathbb R\}.$$

Of course using the multiplicative notation $\chi_\alpha:= e^{\alpha}$  we have that the numbers $c_\beta$  are the coefficients of the generating series:
\begin{equation}\label{gene}
f_B:=\sum_\beta c_\beta e^\beta=\prod_{\alpha\in B}\frac{1}{1-e^\alpha}=\prod_{\alpha\in B}\frac{1}{1-\chi_\alpha}
\end{equation}

Thus the problem is to understand an {\it inversion formula}, the formula  which computes these numbers $c_\beta$ from $f_B$.\smallskip

This is a special case, when $\Delta=\{(1,e^\alpha)\},\ \alpha\in B$, of our theory. 
We have seen, for a general $\Delta$, that, for every function $ f\in Q$  these  numbers are computed as residues by formula (\ref{residuo}).\smallskip

 Our next goal is to give, for any function $$f:={1\over\prod_{(a,\chi)\in\Delta}(1-a \chi )^{h_{a,\chi}}},\ h_{a,\chi}\in\mathbb N$$ and a character $\beta\in C(\Delta)$  a computable expression of $res(\beta^{-1}f)$ in terms of a set (depending on $\beta$) of local residues at the points of the arrangement.\smallskip
This more general setting  gives a weighted version of the partition function which is essentially a variation of the so called Euler--MacLaurin sums of \cite{BV},\cite{BV1}.
 
\subsection{General partial fractions}
 
The results we present from now on are just simple reformulations of the results by Szenes--Vergne \cite{SV1}. We only stress more the role of no--broken circuits which we feel may improve the algorithmic
nature of these formulas.\smallskip

As before, let us take a finite subset $\Delta\subset\mathbb C^*\times\Gamma$. We make two assumptions on $\Delta$:
\begin{enumerate}
\item The   projection $\pi:\Delta\to \Gamma$ gives pairwise non proportional vectors.
\item The characters $\chi$  with as $(a,\chi)$ varies in $\Delta$, span $\Gamma$.
\end{enumerate}
Consider now a pair $(b,\xi)$ with $b\in \mathbb C^*$ and $\xi\in \Gamma_{\mathbb Q}=\Gamma\otimes\mathbb Q$.

We say that $(b,\xi)$ is compatible with $\Delta$ if there is a positive integer $k$ such that $(b^k,\xi^k)\in\Delta$. Notice that by assumption (1), we have that the pair $(b^k,\xi^k)$ is uniquely determined by the pair $(b,\xi)$ and will be called the element of $\Delta$ corresponding to $(b,\xi)$.

We now fix once and for all a total ordering on ${\Delta}$. This ordering induces in an obvious way a partial ordering on the set of pairs compatible with $\Delta$.
We shall consider sequences, possibly with repetitions,
$(b_1,\xi_1),\ldots ,(b_r,\xi_r)$, of elements compatible with $\Delta$ and we shall always assume that
 for each $i=1,\ldots r-1$, either  $(b_i,\chi_i)\leq (b_{i+1},\chi_{i+1})$ or $(b_i,\chi_i)$ and $(b_{i+1},\chi_{i+1})$ are not compatible.
 \begin{definition} A a non broken circuit is:

 a sequence $(b_1,\xi_1),\ldots ,(b_r,\xi_r)$, of distinct elements compatible with $\Delta$, satisfying the following two conditions:
 \begin{enumerate}
\item For each $i=1,\ldots ,r$ there is not a pair $(a,\chi)\in\Delta$ with $(a,\chi)<(a_i,\chi_i)$ and a sequence of non zero integers $m,n_i,\ldots ,n_{r}$ such that:
$$(a\chi)^{m}(b_i\xi_i)^{n_i}\cdots (b_r\xi_r)^{n_r}=1.$$
\item $\xi_1,\ldots ,\xi_r$ are linearly independent.
\end{enumerate}
\end{definition}

We start with another basic identity similar to the one given in 5.1.

\begin{lemma}\label{somma}
\begin{equation}
1-\prod_{i=1}^r z_i=\sum_{\emptyset\subsetneq I\subset \{1,\ldots ,r\}}\prod_{i\in I}(-1)^{|I|+1}(1-z_i)
\end{equation}
\end{lemma}
\proof The proof is by induction on $r$. The case $r=1$ is clear. In general, using the inductive hypothesis, we have
\begin{displaymath}
\sum_{\emptyset\subsetneq I\subset \{1,\ldots ,r\}}\prod_{i\in I}(-1)^{|I|+1}(1-z_i)=
\sum_{\emptyset\subsetneq I\subset \{1,\ldots ,r-1\}}\prod_{i\in I}(-1)^{|I|+1}(1-z_i)-
\end{displaymath}
\begin{displaymath}
 \sum_{I\subset \{1,\ldots ,r-1\}}\prod_{i\in I}(-1)^{|I|+1}(1-z_i)(1-z_r) =
1-\prod_{i=1}^{r-1} z_i \\ + (1-(1-\prod_{i=1}^{r-1}z_i ))(1-z_r)=\end{displaymath}
\begin{displaymath}
=1-\prod_{i=1}^r z_i
\end{displaymath}
\qed

A variation of this formula is the following:
\begin{lemma}\label{somma1} Set $t= \prod_{i=1}^h z_i\prod_{i=h+1}^rz_i^{-1}$, $a\in\mathbb C^*$
\begin{displaymath}
1-t=\sum_{\emptyset\subsetneq I\subset \{1,\ldots ,h\}}\prod_{i\in I}(-1)^{|I|+1}(1-z_i)-a^{-1}\sum_{\emptyset\subsetneq I\subset \{h+1,\ldots ,r\}}\prod_{i\in I}(-1)^{|I|+1}(1-z_i)+\end{displaymath}\begin{equation}\label {cavol}
a^{-1}\sum_{\emptyset\subsetneq I\subset \{h+1,\ldots ,r\}}\prod_{i\in I}(-1)^{|I|+1}(1-z_i)(1-at)
\end{equation}
\end{lemma}
\proof This is immediate from the previous Lemma once we remark that
$ 1-t =1-\prod_{i=1}^hz_i-t(1-\prod_{i=h+1}^rz_i$).\qed

 From this  we get,
\begin{lemma}\label{quozie} Set $t= \prod_{i=1}^h z_i\prod_{i=h+1}^rz_i^{-1}$ with $0\leq h\leq r$. Then
$${1\over \prod_{i=1}^r(1-z_i)}=\sum_{\emptyset\subsetneq I\subset \{1,\ldots ,h\}}{(-1)^{|I|+1}\over (1-t)\prod_{i\notin I}(1-z_i)}-$$\begin{equation}\label{quoz}-\sum_{\emptyset\subsetneq I\subset \{h+1,\ldots ,r\}}({(-1)^{|I|+1}\over\prod_{i\notin I}(1-z_i)}-{(-1)^{|I|+1}\over (1-t)\prod_{i\notin I}(1-z_i)}).\end{equation}
If $a\in\mathbb C^*$ and $ a\neq 1$
$${1\over \prod_{i=1}^r(1-z_i)(1-at)}={a\over a-1}(\sum_{\emptyset\subsetneq I\subset \{1,\ldots ,h\}}{(-1)^{|I|+1}\over (1-at)\prod_{i\notin I}(1-z_i)})-$$\begin{equation}\label{quoz2}-{1\over a-1}(\sum_{\emptyset\subsetneq I\subset \{h+1,\ldots ,r\}}({(-1)^{|I|+1}\over\prod_{i\notin I}(1-z_i)}-{(-1)^{|I|+1}\over (1-at)\prod_{i\notin I}(1-z_i)})+{1\over \prod_{i=1}^r(1-z_i)}).\end{equation}
\end{lemma}
\proof The first relation follows from (\ref{cavol}) dividing by  $(1-t)(1-z_1)(1-z_2)\cdots (1-z_r)$ and taking $a=1$.

The second writing
$$1={a\over a-1}(1-t)-{1\over a-1}(1-at)$$
and then dividing by  $(1-at)(1-z_1)(1-z_2)\cdots (1-z_r)$.\qed

Let us now consider the group algebra $\mathbb C[\Gamma_{\mathbb Q}]$ of the vector space
$\Gamma_{\mathbb Q}$, set $d=\prod_{(a,\chi)\in\Delta}(1-a\chi)$ and $R=\mathbb C[\Gamma_{\mathbb Q}][1/d]$.
\begin{lemma}\label{rela4}
\begin{equation}{n\over 1-x^n}=\sum_{i=0}^n{1\over 1-\zeta^ix},\quad \zeta:=e^{2\pi i/n}.\end{equation}
\end{lemma}
\proof  Take a variable $t$ and derivatives with respect to $t$
$$ {nt^{n-1}\over t^n-x^n}dt=d\log(t^n-x^n)=d\log(\prod_{i=0}^{n-1}(t-\zeta^ix))$$ $$=\sum_{i=0}^{n-1}d\log( t-\zeta^ix))=\sum_{i=0}^{n-1}{1\over ( t-\zeta^ix)}dt.$$ Now set $t=1$ in the coefficients of $dt$.
\qed

 \begin{proposition}\label{expand}  Let  $$\mathcal S:=\{(b_1,\xi_1),\ldots ,(b_M,\xi_M)\}$$ be a sequence of elements compatible with $\Delta$.

Let, for each $j=1,\dots ,M$,  $\zeta_j\in \Gamma_{\mathbb Q}$ be such  that there exists $0\leq n_j\leq m_j$, $m_j>0$, with $\zeta_j^{m_j}=\xi_j^{n_j}$ and set $\zeta=\zeta_1\cdots \zeta_M$.

In $R$ we can write the element $${\zeta \over\prod_{i=1}^M(1-b_i\xi_i)}$$ as a linear combination with constant coefficients of elements of the form $${1\over (1-c_1\psi_1)^{h_1}\cdots (1-c_r\psi_r)^{h_r}}$$
with $h_1,\ldots ,h_r\geq 0$ and
$\{(c_1,\psi_1), \ldots (c_r,\psi_r)\}$ a non broken circuit.
\end{proposition}
\proof Set $D_{\mathcal S}:= \prod_{i=1}^M(1-b_i\xi_i)$.

We first reduce to the case $\zeta=1$. Indeed, for each $j$ take $\theta_j\in \Gamma_{\mathbb Q}$ with $\theta_j^{m_j}=\xi_j$. We then have that $\zeta_j=\theta_j^{n_j}$. Also take a $m_j$-th root of $b_j$, $f_j$ and write
\begin{equation}\label{produ}
1-b_j\xi_j=1-(f_j\theta_j)^{m_j}=\prod_{s=0}^{m_j-1}(1-\exp({2\pi is/m_j})f_j\theta_j)
\end{equation}
Notice that for each $s$, the pair $(\exp({2\pi is/m_j})f_j,\theta_j)$ is compatible with $\Delta$. We now substitute in our sequence $\mathcal S$ the pair $(b_j,\xi_j)$ with the sequence
$\{(\exp({2\pi is/m_j})f_j,\theta_j)\}$, $s=0,\ldots ,m_j-1$. We obtain a new sequence
$$\mathcal S'=\{(d_1,\mu_1),\ldots ,(d_{N},\mu_N)\}$$ with $N=m_1+m_2
+\cdots +m_M$ elements and with $D_{\mathcal S'}=D_{\mathcal S}$.  In this equality, we have replaced $(c_j,\xi_j)$ with a sequence of $m_j$ pairs each having as second coordinate $\theta_j$.
 
 Since
$\zeta_j=\theta_j^{n_j}$ and $n_j\leq m_j$, we have that $\zeta=\mu_1^{\varepsilon_1}\cdots \mu_N^{\varepsilon_ N}$ with $\varepsilon_i\in\{0,1\}$ for each $i$.

If $\zeta\neq 1$, let $j_{\zeta}$ be the least integer such that $\varepsilon_{j_{\zeta}}=1$.
Clearly either $\zeta':=\mu_{j_{\zeta}}^{-1}\zeta$ equals 1 or ${j_{\zeta'}}>{j_{\zeta}}$.

Adding and subtracting $(d_{j_{\zeta}}\mu_{j_{\zeta}})^{-1}\zeta$, a simple computation shows that
$${\zeta\over D_{\mathcal S}}=-d_{j_{\zeta}}^{-1}{\zeta'\over D_{\mathcal S''}}+
d_{j_{\zeta}}^{-1}{\zeta'\over D_{\mathcal S}}$$
where $\mathcal S''=\mathcal S'-\{(d_{j_{\zeta}},\mu_{j_{\zeta}})\}$. If $\zeta'=1$ we are done. Otherwise, everything follows by induction on the number of indices $j$ for which $\varepsilon_j=0$.

Having reduced to the case $\zeta=1$, let us now prove our claim for $1/D_{\mathcal S}$.
Set  Supp($\mathcal S$) equal to the subset in $\mathbb C^*\times \Gamma_{\mathbb Q}$ consisting of pairs $(c,\xi)$ such that $(c,\xi)=(c_j,\xi_j)$ for some $j=1,\ldots ,M$. If Supp($\mathcal S$) is a non broken circuit there is nothing to prove.

Assume  that the first condition in the definition of a no broken circuit is not verified i.e.
there is a pair $(a,\chi)\in\Delta$ and distinct elements $(b_{i_1},\xi_{i_1})<\cdots < (b_{i_t},\xi_{i_t})$ in $\mathcal S$ with $(a,\chi)<(b_{i_1},\xi_{i_1})$ and
 \begin{equation}(a\chi)^{m}(b_{i_1}\xi_{i_1})^{n_{1}}\cdots (b_{i_t}\xi_{i_t})^{n_{t}}=1\end{equation}
for suitable non zero integers $m,n_1,\ldots n_t$. In particular we get, since the elements $\chi$ and $\xi_j$ are characters
 \begin{equation}\label{car}\chi^{m}\xi_{i_1}^{n_{1}}\cdots \xi_{i_t}^{n_{t}}=1\end{equation}
Set $p=|mn_1\ldots n_t|$, and  take $\theta_0,\theta_1,\ldots ,\theta_t$ in $\mathbb C[\Gamma_{\mathbb Q}]$ with the property that $\theta_0^{p/|m|}=\chi$, $\theta_h^{p/|n_h|}=\xi_{i_h}$ for each $h=1,\ldots ,t$. Relation (\ref{car}) becomes
 \begin{equation}\label{cara1}(\theta_0^{\varepsilon_0}\cdots \theta_t^{\varepsilon_t})^p=1.\end{equation}
 with $\varepsilon_j\in\{1,-1\}$. Hence
  \begin{equation}\label{cara2}\theta_0^{\varepsilon_0}\cdots \theta_t^{\varepsilon_t}=1.\end{equation}

Apply, for each $1\leq s\leq t$, the formula \ref{rela4}, for $x=c_{s}\theta_{s}$, $n=p/|n_s|$. Here $c_{s}$ is a
$n$-th root of $b_{i_s}$ and substitute it in $1/D_{\mathcal S}$. We get an expression of  $1/D_{\mathcal S}$ as a linear combination of $p^{t}$ terms each of the form $1/D_{\mathcal S'}$ where $\mathcal S'$ is obtained from $\mathcal S$ substituting each pair $(b_{i_s},\xi_{i_s})$ with  a pair $(c_s,\theta_s)$,
where as above $c_s$ is  a $p/|n_s|$-th root of $b_{i_s}$. In particular ${\mathcal S'}$ has the same cardinality of ${\mathcal S}$ and the sequence of elements in $\Delta$ corresponding to the elements in ${\mathcal S'}$ and ${\mathcal S}$ coincide. Fix such a ${\mathcal S'}$ and take $c_0$ equal to the
$p/|m|$-th root of $a$ such that
 \begin{equation}\label{cara3}(c_0\theta_0)^{\varepsilon_0}\cdots (c_t\theta_t)^{\varepsilon_t}=1.\end{equation}
 We can then apply formula (\ref{quoz}) of Lemma \ref{quozie} and express $1/D_{\mathcal S'}$ as a linear combination of elements of the form $1/D_{\mathcal S''}$ where either the cardinality of $\mathcal S''$ is strictly smaller than that of $\mathcal S'$, or it is the same but $\mathcal S''$ is obtained from $\mathcal S'$
removing a pair $(c_{s},\theta_{s})$ and inserting the smaller pair $(c_0,\theta_0)$.

Assume that the second condition  in the definition of a no broken circuit is not verified i.e.
there are distinct elements $(b_{i_0},\xi_{i_0})<\cdots < (b_{i_t},\xi_{i_t})$ in $\mathcal S$ with $\xi_{i_1},\ldots ,\xi_{i_t}$ linearly dependent. i.e.
 \begin{equation}\label{caru}\xi_{i_0}^{n_{0}}\cdots \xi_{i_t}^{n_{t}}=1\end{equation}
 for suitable integers $n_0,\ldots ,n_t$. Also we can assume that
  \begin{equation}\label{cura} b_{i_0}^{n_{0}}\cdots b_{i_t}^{n_{t}}\neq 1\end{equation}
  otherwise we can use the previous discussion with $(a,\chi)=(b_{i_0},\xi_{i_0})$. As before,  set $p=|n_0\ldots n_t|$, and  take $\theta_0,\ldots ,\theta_t$ in $\mathbb C[\Gamma_{\mathbb Q}]$ with the property that , $\theta_h^{p/|n_h|}=\xi_{i_h}$ for each $h=0,\ldots ,t$. Relation (\ref{caru}) implies
 \begin{equation}\label{carat}\theta_0^{\varepsilon_0}\cdots \theta_t^{\varepsilon_t}=1.\end{equation}
  Applying as before Lemma (\ref {rela4}),    for $x=c_{s}\theta_{s}$, $n=p/|n_s|$. Again $c_{s}$ being a
$n$-th root of $b_{i_s}$, for each $0\leq s\leq t$.  Substituting    in $1/D_{\mathcal S}$, we get an expression of  $1/D_{\mathcal S}$ as a linear combination of $p^{t+1}$ terms each of the form $1/D_{\mathcal S'}$ where $\mathcal S'$ is obtained from $\mathcal S$ substituting each pair $(b_{i_s},\xi_{i_s})$ with  a pair $(c_s,\theta_s)$,
with  $c_s$   a $p/|n_s|$-th root of $b_{i_s}$.

 From relation (\ref{cura}) we deduce that
\begin{equation}\label{cara6}(c_0\theta_0)^{\varepsilon_0}\cdots (c_t\theta_t)^{\varepsilon_t}=a\neq 1.\end{equation}
$a\in\mathbb C^*$.  Now we can apply formula (\ref{quoz2}) in Lemma \ref{quozie} (with $t=(c_0\theta_0)^{-\varepsilon_0}$) and express $1/D_{\mathcal S'}$ as a linear combination of elements of the form $1/D_{\mathcal S''}$ where  the cardinality of $\mathcal S''$ is strictly smaller than that of $\mathcal S'$.
A simple induction then gives our Proposition.\qed

\subsection{The counting formula }  We complete the argument following very closely the proof  given in  \cite{SV1}.
\smallskip

  The setting is thus the following, we have a linear    form $\phi$ with $\langle \phi, \alpha_i\rangle>0,\ \forall i$. Setting $v_i:= {\alpha_i\over \langle \phi, \alpha_i\rangle}$, the vectors $\Psi:=\{v_i\}$ span the  $r-$dimensional real vector space $V$
and   lie  in the  affine hyperplane $\Pi$ of equation $\langle \phi,x\rangle=1$.

 The intersection of  the cone $C(\Psi)=C(\Delta)$ with $\Pi$ is the convex polytope $\Sigma$ envelop of the vectors $v_i$. Each  cone, generated by $k+1$ independent vectors in $\Psi$ (or of $\Delta$),  intersects $\Pi$ in  a $k$ dimensional simplex.

  We have a  configuration of cones   obtained by projecting a configuration of simplices and there is a simple dictionary to express properties of cones in terms of simplices and conversely.

It is well known   that $\Sigma$ is the union of the simplices with vertices independent vectors of $\Psi$.    It is natural to define {\it regular} a point in $\Sigma$  which is not contained in any $r-2$ dimensional simplex (or in the corresponding cone).   The connected components   of the set of regular points are called in \cite{BV1} the {\it big cells}. We complete the set of cells with their boundary cells getting a stratification of $\Sigma$ into cells and a stratification  $\mathfrak S$ of $C(\Delta)$ into polyhedral cones.

 Recall in fact that by  \cite{dp2}, Theorem 5.5 we have .

 \begin{proposition}\label{strat}  Two elements of  $C(\Delta)$ are in the same    cone of $\mathfrak S$   if and only if they are contained in the same set of simplicial cones generated by no broken circuit bases.\end{proposition}
{\bf Remark}\quad  The interesting feature of this statement is that, while $\mathfrak S$ is intrinsically defined, the  no broken circuit bases depend on a choice of a total ordering of $\Delta$.

We fix once and for all an orientation of the vector space  $\Lambda\otimes_{\mathbb Z}\mathbb Q$ and take the invariant $r$-form $\omega_T$ defined in (\ref{can}). 
\begin{definition} Given a chamber $\mathfrak c$ we define the Jeffrey--Kirwan residue, relative to  $\mathfrak c$, of a function $f\in R_T$  to be:
 $$JK(\mathfrak c,f)=(-1)^r\sum_P \sum_S \epsilon_S res_{S,P}(f\omega_T),$$
were $P$ runs over the points of the arrangement and $S$ on the no broken circuit bases in $\Delta_P$ such that \  $\mathfrak c\subset C(S)$ and $\epsilon_S$ equals 1, if $S$ is equioriented with respect to our choice of orientation, -1 otherwise.
\end{definition}\label{JK} 

\begin{remarks} 1) Observe that $JK(\mathfrak c,f)$ does not depend on the choice of the orientation.

2) It follows from  Proposition \ref{periodic}  that $JK(\mathfrak c,\chi^{-1}f)$ as function of $\chi$ is a periodic polynomial.
\end{remarks}
\begin{lemma}\label{independ} Let $\pi: U\to T$ be a finite covering on $T$ of degree $m$. Let $f\in R_T$ then
$f\circ\pi\in R_U$, and
\begin{equation}\label{ugual}JK(\mathfrak c,f)=JK(\mathfrak c,f\circ\pi)
\end{equation}
\end{lemma}
\proof A non broken circuit  $S$ associated to a point $P$  in $T$ is also associated to each point in $\pi^{-1}(p)$ in $U$. In this way one obtains all points in $U$ associated to $S$.
Also from the definition of the local residue and Lemma \ref{multipl}, we deduce that if $Q\in \pi^{-1}(P)$
$$res_{S,Q}(f\omega_U)=m^{-1}res_{S,P}(f\omega_T).$$
Since $m=|\pi^{-1}(P)|$ everything follows.\qed

 \begin{lemma} 1) Consider  a  non broken circuit  basis  $S:=\{\psi_1,\dots,\psi_r\}$, a character $\mu=\prod_i\psi_i^{k_i}$ with $0<k_i\in \mathbb N$ and a function  $$\gamma:={\mu^{-1}\over \prod_{i=1}^k(1-a_{i_1}\psi_{i_k})^{h_i}}, \ 1\leq i_1<\dots<i_k\leq r.$$
 Then 
 \begin{equation}\label{finn}res(\gamma)=(-1)^r\epsilon_S\sum_{P}res_{S,P}(\gamma\omega),\end{equation}where $P$ lies in the finite set of points associated to $S$ (i.e. defined by the equations $1-a_i\psi_i=0$, $i=1,\ldots ,r$). In particular if $k<r$ both terms are zero.
 
 2) If $\mu$ is regular and $\mathfrak c$ is the unique chamber containing $\mu$, $$(-1)^r\epsilon_S\sum_{P}res_{S,P}(\gamma\omega)=JK(\mathfrak c,\gamma).$$
The sum being on the set of $P$ associated to $S$. \end{lemma}
 \proof
 1) Taking as $\Delta$ the set $\Delta_S=\{(a_1,\psi_1),\dots,(a_r,\psi_r)\}$ we see that
 $$JK( C(S),\gamma)=(-1)^r\epsilon_S\sum_{P}res_{S,P}(\gamma\omega),$$
so by Lemma \ref{independ} the right hand side of    formula (\ref{finn}) is independent of a torus on which the $\psi_i$  are defined. By definition the same is true for $res(\gamma)$. Thus we can assume we are on
the torus for  which $S$ is a basis of the character group. In this case, there is a unique point $P$ associated to $S$. Furthermore  $\omega=\epsilon_Sd\log\psi_1\wedge\dots\wedge d\log\psi_r$.   Now one separates the variables and reduces to the 1 dimensional   case. This    is an elementary immediate instance  of the knapsack problem (see \S  1.1 or \cite{SV1}).
 
 2) We need to see that, if $(S',Q)$ is not one of our $(S,P)$, then 
 $res_{S',Q}(\gamma\omega)=0.$ If $Q\neq P$ at least one on the factors $1-a_i\psi_i$  is holomorphic in $Q$ and we can apply part 1) of Theorem \ref{primc}. If $P=Q$   remark  that the form $\gamma\omega$ is in the image of the cohomology of $\mathcal A_{\Delta_S}\supset \mathcal A_\Delta$. Our claim now follows from formula (\ref{clocal}) of Theorem \ref{loc} applied to $P\in \mathcal A_{\Delta_S}$.
 \qed

\begin{theorem}\label{regcase}  Let $$f:={1\over\prod_{(a,\chi)\in\Delta}(1-a \chi )^{h_{a,\chi}}},\ h_{a,\chi}\in\mathbb N$$
If $\beta$ is in the closure  of  a chamber $\mathfrak c$:
\begin{equation}\label{gfor}c_\beta=res(\beta^{-1} f \omega)=JK(\mathfrak c,\beta^{-1} f  )\end{equation}
\end{theorem}
\proof  Using Lemma \ref{independ}, we can pass to any finite covering $U$ of $T$. For a suitable covering we can find a character $\xi$ such that $\beta\xi$ lies in the interior of $\mathfrak c$ and the function 
$\xi f$ has an expansion as in Proposition \ref{expand}. Having done this, we are reduced to prove our identity for each single term of this expansion. Write $\beta^{-1}f=(\beta\xi)^{-1}\xi f$.
Each term of the expansion of  $\xi f$ satisfies the hypotheses of the previous lemma with $\mu=\beta\xi$. The result follows.\qed

Remarks:  From the remark after Definition \ref{JK} follows that in the previous setting, $c_\beta$ as  function of $\beta$ is a periodic polynomial. It has also a remarkable continuity property.  That is if $\beta$ is in the boundary of two different chambers  there are two possibly different periodic polynomials on these chambers which agree on the intersection of their closure, in particular on $\beta$.

\end{document}